\documentclass[11pt, a4paper]{amsart}
\input xy   
\xyoption{all}
\numberwithin{equation}{section}
\swapnumbers

\theoremstyle{plain}
  \begingroup
        \newtheorem{theorem}[equation]{Theorem}
        \newtheorem{lemma}[equation]{Lemma}
        \newtheorem{proposition}[equation]{Proposition}
        \newtheorem{corollary}[equation]{Corollary}
        \newtheorem{fact}[equation]{Fact}

        \newtheorem{remark}[equation]{Remark}
	\newtheorem{definition}[equation]{Definition}
        \newtheorem{notation}[equation]{Notation}
        \newtheorem{convention}[equation]{Convention}
        \newtheorem{sinnadaitalica}[equation]{}
  \endgroup

\theoremstyle{definition}
  \begingroup
        \newtheorem{example}[equation]{Example}
        \newtheorem{comment}[equation]{Comment}
        \newtheorem{sinnadastandard}[equation]{}
        \newtheorem{observation}[equation]{Observation}
  \endgroup

\newcommand{\mr}[1]{\buildrel {#1} \over \longrightarrow}
\newcommand{\ml}[1]{\buildrel {#1} \over \longleftarrow}

\newcommand{\colim}[2] {\displaystyle \lim_{\overrightarrow{#1}} {#2}}

\newcommand{\rimply}{\Rightarrow}

\newcommand{\cc}{\mathcal}
\newcommand{\bb}{\mathbb}
\newcommand{\ff}{\mathsf}

\newcommand{\mmr}[1]{\buildrel {#1} \over \hookrightarrow}

\begin{document}

% $$\colim{U \in \cc{F}_p}{X_U}$$ 
% $$\colimite{U \in \cc{F}_p}{X_U}$$
\title{Representation theory of mv-algebras}

\author{Eduardo J. Dubuc and Yuri A. Poveda}

\begin{abstract}
In this paper we develop a general representation theory for
mv-algebras. We furnish the appropriate categorical background to study
this problem. Our guide line is the theory of classifying topoi of
coherent extensions of universal algebra theories. Our main result
corresponds, in the case of mv-algebras and mv-chains, to the
representation of commutative rings with unit as rings of global
sections of sheaves of local rings. \emph{We prove that any mv-algebra is
isomorphic to the mv-algebra of all global sections of a sheaf of
mv-chains on a compact topological space}. This result is intimately
related to McNaughton's theorem, and we explain why our representation
theorem can be viewed as a vast generalization of McNaughton's. On
spite of the language utilized in this abstract, we wrote this paper in a way that, we hope, could be read without much acquaintance with either sheaf theory or mv-algebra theory.  
\end{abstract}

\maketitle

{\sc preface}

We wrote this paper in a way that, we hope, could be read without
much acquaintance with either sheaf theory or mv-algebra theory. Our basic reference on mv-algebras is  the book  'Algebraic Foundations of Many-valued Reasoning' \cite{COM}, and we refer to this book and not to the original sources for the known specific results we utilize. Of sheaf theory we need only some basic general facts, the reader may consult the book  'Sheaves in Geometry and
  Logic' \cite{MM}. The basic reference for category theory is of course the classical  'Categories for the working mathematician' \cite{M}.

\tableofcontents

{\sc introduction.} 

In this paper we develop a general representation theory for mv-algebras as algebras of global sections of sheaves of mv-chains. Our work generalizes previous results in this direction, like the representation theorem for locally finite mv-algebras proved in \cite{CDM}. Here we furnish the appropriate categorical background to study the representation theory of mv-algebras. Our guide line is the theory of classifying topoi first developed in the case of rings by M. Hakim, and thereafter placed in the general context of universal algebra by category theorists, in particular M. Coste  \cite{C}. Our main result corresponds, in the case of mv-algebras and mv-chains, to the representation of commutative rings with unit as rings of global sections of sheaves of local rings, as in, for example,  R. Hartshorne  \cite{H}. We prove that \emph{every mv-algebra is isomorphic to the mv-algebra of all global sections of its prime spectrum} (theorem \ref{main}). We analyze and develop carefully the various steps that lead to this theorem.

In section \ref{terminology} we recall the basic fundamental results on sheaves and on mv-algebras, and in this way fix notation, terminology,  and set up a place for reference.  

In sections  \ref{coZariski} and \ref{duality} we develop the basic framework of the representation and duality theory for mv-algebras.  
The localization of a ring in a prime ideal is a covariant
construction in the sense that if we have two prime ideals $P \subset
Q$, we have morphism of local rings $A_P \to A_Q$. In the case of
mv-algebras, the quotient of an mv-algebra by a prime ideal is a
contravariant construction, in the sense that we have a morphism of
mv-chains in the other direction $A_{/Q} \to A/_P$. This indicated to
us that the appropriate topology for the set of prime ideals in the
case of mv-algebras is not the Zariski topology  as in the case of
rings, but a topology that we call \emph{coZariski}. Its base of open
sets is given by, for each $a \in A$, the sets $\{P\,|\, a \in P\}$, and not by the sets $\{P\,|\, a \notin P\}$ as in the Zariski topology. We think that the consideration of the Zariski topology in the representation theory of mv-algebras  has blocked the development of the theory beyond the particular case of hyperarchimedean algebras, where the two topologies coincide. We construct the prime spectrum sheaf of a mv-algebra following standard methods of sheaf theory, and introduce the notion of \emph{mv-spaces} and its corresponding category mimicking the algebraic geometry notion of ringed space. Our main result (theorem \ref{main}) means that \emph{the category of mv-algebras is the dual of a category of mv-spaces.} This is similar to Grothendieck method in algebraic geometry, which is based in the consideration of the dual of the category of rings as the category of affine schemes. 
We show that the main theorem follows from two facts. First, that the
prime spectrum is a compact topological space, the \emph{compacity
  lemma} (lemma \ref{compacity}), and second, what we call \emph{pushout-pullback lemma} (lemma \ref{pushpull}). This lemma means that given two elements $a_1, \,a_2$ in a mv-algebra $A$, any two other elements $b_1, \,b_2$ such that $b_1 = b_2 \; mod(a_1 \vee a_2)$, can be "glued" into a single element $b$, unique $mod(a_1 \wedge b_2)$, such that $b = b_1\;mod(a_1)$, and $b = b_2\;mod(a_2)$. This lemma, not surprisingly, is intimately ligated with McNaughton's theorem. 

In section \ref{msofamvalgebra} we consider the maximal spectrum
furnished with the coZariski topology, and the $[0, \,1]$-valued
morphism spectrum furnished with the Zariski topology, as studied in
\cite{CDM}. These two spectral spaces have the same underlying set,
but the coZariski topology is always finer than the Zariski
topology, and strictly so, unless the algebra is hyperarchimedean. We
prove several results necessary in the proof of the pushout-pullback
lemma, and and helpful to study the relationship of this lemma with \mbox{McNaughton} theorem. We are lead to a concept of strong semisimplicity \mbox{(definition \ref{ssemisimple}).} Strongly semisimple mv-algebras are intermediate between hyperarchimedean algebras and general semisimple algebras. The primarily example are the finitely presented mv-algebras. 

In section \ref{hyperarchimedean}, 
before attacking the general case,  we  consider the case of
hyperarchimedean algebras, and prove the main representation theorem for these algebras. Here every prime ideal is
maximal, the Zariski and coZariski topologies coincide, and the prime spectrum  is separated (Hausdorff).
Furthermore, all the fibers in the prime spectrum are subalgebras of the 
real unit interval $[0,\, 1]$. All these facts allow us to apply classical
methods as in \cite{CDM}. We generalize to arbitrary hyperarchimedean algebras the duality obtained in \cite{CDM} for locally finite mv-algebras.

In section \ref{fp} we prove and/or recall several results on finitely presented mv-algebras, and prove the pushout-pullback lemma for these algebras. A key result is the gluing lemma 5.3 of \cite{NPM}, that we adapt and prove in our context in lemma \ref{pushpullfp}.

In section \ref{pushpulls} we prove the general case of the pushout-pullback lemma. It follows from the case of finitely presented algebras by categorical nonsense because finite limits commute with filtered colimits. We do here an sketch of an explicit proof in the case of pullbacks of mv-algebras. 

In section \ref{posetcoZariski} we prove the compacity lemma, that is, we prove that the prime spectrum furnished with the coZariski topology is a compact topological space. We develop (following general guide lines coming from the theory of Grothendiek topologies)  a construction of its lattice of open sets without constructing the underlying set first. This yields a compact locale.  Then the prime spectrum is the space of points of this locale. Such space will be be a compact topological space provided the locale has enough points, a fact that follows by an standard application of Zorn's Lemma.

In section \ref{McN} we deal with McNaughton's theorem. We show
explicitly, and prove, that this theorem is equivalent to our
representation theorem for the particular case of free mv-algebras. We
show that a finite open cover of the prime spectrum of the free
mv-algebra on $n$ generators yields (by restricting  the open sets of
the cover to the maximal spectrum), a finite decomposition of the $n$
dimensional cube by convex polyhedra. Once this is understood, an
isomorphism between the mv-algebra of McNaughton functions and the
mv-algebra of global sections of the prime spectrum of the free algebra
becomes evident. This shows the equivalence of the two results. In
particular, we obtain a proof of McNaughton's theorem as an
application of our representation theorem. These results also furnish
a conceptual context for the theorem. They exhibit the "local" nature
of the concept of McNaughton function as an instance of the usual
topological notion of localness. Our representation theorem can be
seen as a vast generalization of McNaughton's theorem, from free
mv-algebras to arbitrary mv-algebras. 

Finally, in an appendix we prove some general
results as we need them in section \ref{posetcoZariski}. This results conform what we call \emph{Sheaf theory of posets}, they are
known and can be found (in one form 
or another) in the literature. The development we do here is original
in print as far as we know. We did not find any treatment which arrive fast to the results, and that in the statements and the proofs is  
faithful to the idea we want to stress here, namely, that of a reflexion of the theory of sheaves on $Set$-based categories into \mbox{$\{0,\,1\}$-based} categories, namely, partially ordered sets. Here, posets play the role of categories, the two
elements poset play the role of the category of sets,  
\mbox{inf-lattices} play the role of categories with finite limits, and locales
play the role of Grothendieck topoi.

%\vspace{1ex}

{\sc acknowledgments}
The first author is grateful to Roberto Cignoli for several helpful
and inspiring
conversations on the subject of this article, and the second author is
grateful to Manuela Busaniche for the important help she gave him on
several specific questions on mv-algebras.

%\tableofcontents

%%%%%%%%%%%%%%%%%%%%%%%%%%%%%%%%%%%%%%%%%%%%%%%%%%%%%%%%%%%%%%%%%%%%%%
\section{Background, terminology and notation}
\label{terminology} 
%%%%%%%%%%%%%%%%%%%%%%%%%%%%%%%%%%%%%%%%%%%%%%%%%%%%%%%%%%%%%%%%%%%%%%
In this section we recall some basic mv-algebra and sheaf (on 
topological space) theory, and in this way fix notation and
terminology.

\vspace{1ex}

For the definition and basic facts on mv-algebras the reader is
advised to have at hand the book ``Algebraic foundations of
many-valued reasoning'', reference \cite{COM}.

\begin{sinnadastandard} {\bf MV-Algebras}  \label{mvfacts1} \end{sinnadastandard}

MV-algebras are models of an equational theory in universal algebra.

\vspace{1ex}

\begin{enumerate}
\item 
A mv-algebra $A$ has a $2$-ary, a
$1$-ary  and a $0$-ary primitive operations, denoted $\oplus,\;
\neg,\; 0$, subject
to universal axioms.

\vspace{1ex}

\item 
It is convenient to introduce one $0$-ary and two $2$-ary derived
operations, defined by the following formulae. 
$$
1 = \neg 0,\;\;\; x \odot y = \neg(\neg x \oplus \neg y), \;\;\; x \ominus y = x \odot \neg y
$$ 

\vspace{1ex}

\item \label{nx} Given an integer $n$, we let $nx= x \oplus x \oplus x \cdots
  \oplus x$, $n$ times, and $0x = 0$.

\vspace{1ex}

\item
There is a partial order relation defined by $x \leq y \iff \exists
z,\, x \oplus z = y$

\vspace{1ex}

\item
A useful characterization is the following:$\; x \leq y \iff x \ominus y
= 0$.

\vspace{1ex}
 
\item The partial order is a lattice, with \emph{supremun} denoted $ \vee$ and
  \emph{infimun} denoted $\wedge$. These lattice structure is also
  given by formulae. The bottom and top elements are
  $0$ and $1 = \neg 0$. 

\vspace{1ex}

\item It holds  $\; x
  \wedge y \leq x \oplus y\;$ and $\; x \odot y \leq x \wedge y$.

\vspace{1ex}  

\item \label{distance} 
There is a \emph{distance} operation defined by $d(x,\,y) =
  (x\ominus y) \oplus (y \ominus x)$, and it holds $d(x,\,y) = 0 \;
  \iff \; x = y$

\vspace{1ex}

\item \label{p1}
The following equations holds in any mv-algebra:

$x \oplus \neg x = 1$, $\;(x \ominus y) \wedge (y \ominus x) = 0$

\vspace{1ex}

\item \label{p2}
Given any mv-algebra $A$, $x,\, y \, \in A$ and any integer $n$,

$n(x \wedge y) = nx \wedge ny$. It follows 
$x \wedge y = 0 \;\rimply\; nx \wedge ny = 0$.

%\qed\end{facts}

\vspace{1ex}

\item $A = \{0\}$ is the trivial mv-algebra, $A$ is said
\emph{nontrivial} if $1 \neq 0$.  

\vspace*{3ex}

\item \label{[01]} The closed real unit interval $[0, \, 1] \subset \bb{R}$ is the basic
example of an mv-algebra. The structure is given by:

$x \oplus y = min(1, \, x + y), \; \neg x = 1 - x, \; x \ominus y =
max(0, \, x-y),$ \mbox{$x \odot y = max(0, \, x+y-1)$}, $\;d(x,\, y) =
|x - y|.$
\end{enumerate}
\begin{sinnadastandard} \label{mvfacts2}
{\bf Morphisms of mv-algebras} 
\end{sinnadastandard}
\begin{enumerate}
\item \label{mvalgebras} We shall denote by $mv\cc{A}$ the category of
  mv-algebras. It 
  has all the standard constructions and 
properties of any variety of universal algebras. Morphisms of
mv-algebras are defined as usual, they are functions 
which preserve the three primitive operations. It follows that a
morphism preserves also all the derived operations since
these operations are defined by formulae. In particular it preserves
\mbox{$0$ and $1$.} Thus the zero function $A \mr{0} B$ is a morphism only
when $B$ is the trivial algebra. 

\vspace{1ex}
                                                                                                                      
\item Given a morphism of mv-algebras $A \mr{\varphi} B$, the
\emph{image}, $Im(\varphi) \subset B$, is a subalgebra of $B$,
\mbox {$Im(\varphi) = \{y \in B \; |
\; \exists x \in A, \; \varphi(x) = y \}$}.

\vspace{1ex}

\item Congruences are associated with ideals. 
  
An \emph{ideal} of an mv-algebra $A$ is a subset $I \subset A$
satisfying: 

I1) $0 \in I$

I2) $x \in I,\; y \in A,\; y \leq x,\; \rimply \; y \in I$

I3)  $x \in I, \; y \in I,\; \rimply \; x\oplus y \in I$

An  ideal is said to be  \emph{proper} if $1 \notin I$. 

\vspace{1ex}

\item \label{aplusbequalab}

Given a finite number of elements in an mv-algebra,  we write
$(x,\, y \, 
\ldots \, z)$ for the ideal generated by $x,\, y \, \ldots \, z$. In
any mv-algebra it holds 
$(x, \, y) = (x \oplus y)$. Consequently, all finitely generated
ideals are \emph{principal}.

\vspace{1ex}

\item \label{distance2} Given a morphism $A \mr{\varphi} B$,   the
\emph{kernel}, $Ker(\varphi) \subset A$, is an ideal of A, $Ker(\varphi) 
= \{x \in A \; | \; \varphi(x) = 0\}$. Any ideal $I
\subset A$
determines a congruence by stipulating $x \sim y \; \iff \; d(x,\,y)
\in I$. The quotient is denoted $A \mr{\rho} A/I$. It holds $I =
Ker(\rho)$. Given $x \in A$, we denote $[x]_I \in A/I$, $[x]_I =
\rho(x)$. Thus, $[x]_I = [y]_I  \; \iff \; d(x,\,y) \in I$.

\vspace{1ex} 

\item \label{primetheorem} 
Consider now prime ideals $P$ (see (\ref{prime}) below). The following
  equivalence holds in any mv-algebra (see \cite{COM} 1.2.14): 
\vspace{1ex}

$((\forall P\; prime)\; x \in P  \rimply  y \in P)\; \iff \; (y) \subset (x)$
\end{enumerate}
\begin{sinnadastandard} \label{mvfacts3}
 {\bf MV-chains and Prime ideals}
\end{sinnadastandard}
\begin{enumerate}
\item A mv-algebra is a \emph{mv-chain} if it is non trivial and its
  order relation is 
total. Thus, mv-chains are characterized by the (non universal)
axioms:

C1) $1 \neq 0$.

C2) $x \ominus y = 0 \;\; or \;\; y \ominus x = 0$.

\vspace{1ex}

\item From (\ref{p1}) it follows that a mv-algebra $A$ is a mv-chain if and
  only if it satisfies:

C'1) $A \neq \{0\}$.

C'2) $x \wedge y = 0 \;\; \rimply \;\; x = 0 \; or \; y = 0$

\vspace{1ex}

\item \label{prime} An ideal $P$ is \emph{prime} if it satisfies the following
conditions:

P1) $1 \notin P$.

P2) For each $x, \, y $ in $A$, either $(x \ominus y) \in P$ or $(y
\ominus x) \in P$.

\vspace{1ex}

\item From (\ref{p1}) it follows that an ideal $P$ is prime if and
  only if it satisfies:

P'1) $P \neq A$.

and either of the following two equivalent conditions

P'2) $x \wedge y \in P \;\; \rimply \;\; x \in P \; or \; y \in P$

P'2) $x \wedge y = 0 \;\; \rimply \;\; x \in P \; or \; y \in P$

\vspace{1ex}

\item \label{chainprime} The following then holds by definition:

Given any morphism  $A \mr{\varphi} B$, then $Ker(\varphi)$ is a prime ideal
iff $Im(\varphi)$ is a mv-chain.

\vspace{1ex}

\item \label{phimenosunoprime}

Given any morphism  $A \mr{\varphi} B$ and a prime ideal $P$ of $B$, the ideal $\varphi^{-1}(P)$ is a prime ideal of $A$. 

\end{enumerate}
\begin{sinnadastandard} \label{mvfacts4}
{\bf Simple mv-algebras and Maximal ideals}
 \end{sinnadastandard}
\begin{enumerate}
\item  \label{maximal0}
Maximal ideals are prime.Given a prime ideal $P$, there exist a
maximal ideal $M \supset P$ (see \cite{COM} 1.2.12). 

\vspace{1ex} 

\item \label{maximal}
An ideal $M$ is maximal if and only if 
$M \neq A$ and, for each $a \in A$, $\; a \notin M \; \iff \;
 \exists\, n \geq 1 \;\; \neg\, nx \in M$ (see \cite{COM} 1.2.2).

\vspace{1ex}

\item  \label{maximal2}
Recall that a mv-algebra is \emph{simple}
iff it is nontrivial and 
$\{0\}$ is its only proper ideal. By definition, 
an ideal $M$ is maximal if and only if the quotient
 $A/M$ is a simple mv-algebra.

\vspace{1ex}

\item \label{uniquesubalgebra}
Given any maximal ideal $M$ of a mv-algebra $A$: 

a) $A/M$ is isomorphic to a subalgebra of
the real unit interval $[0,\, 1]$ (see \cite{COM} 1.2.10, 3.5.1). 

\vspace{1ex}

b) There is a unique
embedding $A/M \mmr{} 
[0,\,1]$ which determines by composition a morphism $\chi_M: A \mr{} A/M
\mmr{} [0,\,1]$. Furthermore $ker(\chi_M) = M$, and $\chi_{ker(\chi)}
= \chi$ (see \cite{COM} 7.2.6).

\vspace{1ex}

\item \label{phimenosunomaximal}
Given any morphism  $A \mr{\varphi} B$ and a maximal ideal $M$ of $B$, the ideal $\varphi^{-1}(M)$ is a maximal ideal of $A$ (see \cite{COM} 1.2.16). 

\vspace{1ex} 

Contrary to \ref{mvfacts3} (\ref{phimenosunoprime}), this fact distinguishes mv-algebras from other related structures, and follows from the characterization \ref{mvfacts4} (\ref{maximal}).

\end{enumerate}

\vspace{1ex}

\begin{sinnadastandard} {\bf Presheaves and
    Sheaves} \label{sheavesfacts} \end{sinnadastandard} 
For the following the reader may consult the book ``Sheaves in Geometry and
Logic'', reference \cite{MM}.
 
Given a topological space X, we denote $\cc{O}(X)$ the lattice of open
sets. 
\begin{enumerate} 
\item A \emph{presheaf} $F$ is a contravariant functor $\cc{O}(X)^{op}
\mr{F} \cc{S}$ into the category $\cc{S}$ of sets. Given $U \subset V$
in  $\cc{O}(X)$, and $s \in FV$, we denote $s|_U \in FU$ the image of
$s$ under the map $FV \to FU$.

\vspace{1ex}

\item A \emph{sheaf} is a presheaf $F$ such that
given any open cover $U_\alpha \subset U$ in $\cc{O}(X)$, the
following \emph{sheaf axiom} holds for $F$:

$
\hspace{-5ex}
\forall \, s_\alpha \in FU_\alpha, \; such \; that \;\; s_\alpha |_{U_\alpha
  \cap U_\beta} \;=\; s_\beta |_{U_\alpha \cap U_\beta}, \;\; \exists !
\, s \in FU, \; s|_{U_\alpha} = s_\alpha. 
$

\vspace{1ex}

\item An \emph{etale space} over $X$ is a local homeomorphism $E
  \mr{\pi} X$.
\begin{enumerate}
\item  Given $p \in X$, the \emph{stalk} (or \emph{fiber}) over $p$ is the set
  $E_p = \pi^{-1}(p)$.

\vspace{1ex}

\item Given $U \in \cc{O}(X)$, a \emph{section} defined on $U$ is a continuous
  function $U \mr{\sigma} E$, $\pi \sigma = id$. If $U = X$, $\sigma$
    is said
  to be a \emph{global} section. The set of sections
  defined on $U$ is denoted $\Gamma(U,\,E)$. 
\end{enumerate}

\vspace{1ex} 

\item \label{global} An etal space $E \to X$ is said to be
  \emph{global} if each point $e \in E$ is the image of some global
  section.

\item \label{T_2sheaf} Given any two sections $\sigma$, $\mu$, 
  the set $[\lbrack \sigma = \mu \rbrack] = \{ x \;|\; \sigma(x) = \mu(x)
\}$ is open, and it is also closed (hence
\emph{clopen}) when $E$ is Hausdorff.

\vspace{1ex}

\item The assignment $U \mapsto \Gamma(U,\,E)$ together with the
  restriction maps determine a sheaf $\;\Gamma(-,\,E)\;$ (the sheaf of
  sections of $E \to X$).

\vspace{1ex}

\item \emph{(Godemment construction)} Given any presheaf $\cc{O}(X)^{op} \mr{F}
  \cc{S}$, we can associate with $F$ an etal space $E_F \mr{\pi}
  X$. 

For each $p \in X$, we can take the colimit of the system 
$\{FU\}_{U \in \cc{F}_p}$ indexed by the filter $\cc{F}_p$ of open neighborhoods
  of $p$. Then, set $(E_F)_p = \colim{U \in \cc{F}_p}{FU}$. We have maps $FU
  \to (E_F)_p$. Given $U \ni p$ and $s \in FU$, we denote $[s]_p$ the
  corresponding element in $(E_F)_p$.  Abusing this notation, we have 
$(E_F)_p = \{[s]_p \; | \; s \in FU,\; U \ni p\}$. 

We define the
  set $E_F$ as the disjoint union of the sets $(E_F)_p$, $p \in
  X$.  
$E_F = \{([s]_p,\,p),\; s \in FU,\; U \ni p,\; p \in X  \}$. 
The map $E_F \mr{\pi} X$ is defined to be the projection,  $\pi([a]_p,\,p)
= p$.

Given $U \in \cc{O}(X)$, each element $s \in FU$ determines a section 
$\widehat{s}: U \to E_F$, $\widehat{s}(p) = ([s]_p,\, p)$. We
topologize the set $E_F$ by taking as a base of open sets all the
image sets $\widehat{s}(U)$. Under this topology the map $\pi$ becomes
a local homeomorphism.

\vspace{1ex}

\item \label{basesobra}For the Godemment construction it is equivalent to use any
  \emph{filter base} $\cc{B}_p \subset \cc{F}_p$ in place of
  $\cc{F}_p$ (notice that $\cc{B}_p$ is cofinal in  $\cc{F}_p$). 

\vspace{1ex}

\item The assignment  $s \mapsto \widehat{s}$ defines a natural transformation
  \mbox{$FU \to \Gamma(U,\,E_F)$.} This transformation is an isomorphism iff
  F is a sheaf. 

\vspace{1ex}

\item \label{sheavescat} We shall denote by $\cc{S}h(X)$ the category of sheaves
  on a (fixed) topological space $X$. It is a Topos. The objects are
  etale spaces $E 
  \mr{p} X$, and a morphism from $E \mr{p} X$ to $G \mr{q} X$ is a continuous
  function $E \mr{f} G$ such that $q\circ f = p$. The sheaf of
  sections and the Godemment construction establish an equivalence
  between this category and the category which has as objects the
  Presheaves $\cc{O}(X)^{op} \to \cc{S}et$ which satisfy the sheaf
  axiom, and as morphisms the natural transformations of functors.    
\end{enumerate}

\vspace{1ex}

{\bf Convention:} \emph{When dealing with sheaves we shall use the
contravariant functor datum  or the etale space datum indistinctly,  
and we shall say \emph{``sheaf''} in both cases. When the base space
is not fixed, we shall write $\ff{E} = (X,\, E)$.} 

\vspace{1ex}

Sheaves with variable base space also form a category, but in a
different way (see definition \ref{mv-spaces})

%%%%%%%%%%%%%%%%%%%%%%%%%%%%%%%%%%%%%%%%%%%%%%%%%%%%%%%%%%%%%%%%%%%%%%
\section{The prime spectrum of mv-algebras} \label{coZariski} 
%%%%%%%%%%%%%%%%%%%%%%%%%%%%%%%%%%%%%%%%%%%%%%%%%%%%%%%%%%%%%%%%%%%%%%

{\sc the prime spectrum $\ff{Spec}_A$} \label{sofamvalgebra}  

Given a mv-algebra $A$, we can associate with $A$ a topological space
$Z_A$, and a global sheaf over $Z_A$, $\; \ff{Spec}_A = (E_A \mr{\pi} Z_A)$,
as follows:

\vspace{1ex}

%{\bf The prime spectrum $\ff{Spec}_A$ .} \label{psofamvalgebra}

\begin{sinnadastandard}{ \bf Construction of $\ff{Spec}_A$ .} \label{Spec_A} 
 
The set of points of $Z_A$ is the set of all prime ideals $P \subset A$. For
each \mbox{$a \in A$,} let $W_a \subset Z_A$ be the set $W_a = \{P \;|\; P \ni
a\}$. It is immediate to check  $W_a \cap W_b = W_{a \oplus b}$, so
we can take these sets as a base of a topology, the \emph{coZariski topology}
(notice that $W_0 = Z_A$ and $W_1 =  \emptyset$).

\vspace{1ex}

We define the set $E_A$ as the disjoint union of the mv-chains
$A/P$, $P \in Z_A$. That is $E_A = \{([a]_P,\,P),\; a \in A,\; P \in
Z_A\}$. The map $E_A \mr{\pi} Z_A$ is defined to be the projection,
$\pi([a]_P,\,P) = P$.

Each
element $a \in A$ defines a global 
section (as a function of sets) \mbox{$Z_A \mr{\widehat{a}} E_A$,}
$\widehat{a}(P) = ([a]_P,\, P)$.   
\begin{observation} \label{notation}
By the definitions given we have: $W_a = [[\,\widehat{a} = 0\,]]$ and
$W_{d(a,\,b)} = 
[[\,\widehat{a} = \widehat{b}\,]]$ (see \ref{mvfacts2}
(\ref{distance2}), \ref{sheavesfacts} (\ref{T_2sheaf}))
\qed\end{observation}
\emph{Thus, the open sets basis $W_a$ is formed by the Zero sets of the sections
  $\widehat{a}$.} We shall also use the notation $Z(\widehat{a}) = [[\,\widehat{a} = 0\,]]$.

\vspace{1ex}

 We shall now
define a topology in $E_A$ in such a way that $\pi$ becomes a local
homeomorphism, or \emph{etale space} over $Z_A$.  
The open base for this topology consists of  all the
image sets $\widehat{a}(W_b)$. These sets are closed under
intersections. In fact, given $\widehat{a}(W_b)$ and 
$\,\widehat{c}(W_d)$, 
we have  
%$W_{b \oplus d \oplus d(a,\, c)} \subset W_b \cap W_d$, 
$\widehat{a}(W_b)\cap \widehat{c}(W_d)=
\widehat{s}(W_{b \oplus d \oplus d(a,\, c)})$ (where $s$ is
either $a$ or $b$). 

It is easy, and left to the reader, to check that under this topology
$\pi$ becomes a local homeomorphism and  the global sections $\widehat{a}$
continuous and open functions. Furthermore, this topology is the final
topology with respect to the functions $\widehat{a}$. 
The sections $\widehat{a}$ show that $\ff{Spec}_A$ is a global sheaf.
\qed\end{sinnadastandard}

\begin{sinnadastandard} \label{sections}
{\bf Sections of $\ff{Spec}_A$.}
Given any open set $U \subset Z_A$, the sections over $U$ of
$\ff{Spec}_A$ are the sections of $\pi$ (as functions of sets)  $U
\mr{\sigma} Z_A$  which
locally are of the form $\widehat{a}$. We explain: Given any $P \in
U$, $\sigma{P} = ([s]_P, \, P)$ 
for some $s\in A$. Then, the fact that $\sigma$ is in 
$\Gamma(U,\, E_A)$ means that there is an $a \in A$, $W_a \ni P, W_a \subset U$
  such that 
 $\sigma(Q) = ([s]_Q,\,Q)$, the same $s$ for all $Q \in W_a$. That is,
  $\sigma = \widehat{s}$ on $W_a$. This is equivalent to the existence
  of a compatible family of
sections $\widehat{s_i}$ defined on a open covering $W_{a_i}$ of $U$, such
that $\sigma|_{W_{a_i}} = \widehat{s_i}|_{W_{a_i}}$ for all $\,i\,$.
\end{sinnadastandard}

\vspace{1ex}
  
The construction above is reminiscent of the Godemment
construction. We show now that it actually is equivalent to the
Godemment construction on a certain presheaf.

\begin{sinnadastandard}{ \bf Godemment construction for
    $\ff{Spec}_A$.} \label{spec_Abis}
The topological space $Z_A$ is defined as before. We consider now a
presheaf $\cc{O}(Z_A)^{op} \mr{} \cc{S}$. This presheaf is defined on
the base of open sets by the assignment $W_a \mapsto A/(a)$, where
$(a)$ indicate the principal ideal generated by $a$. From
\ref{mvfacts3} (\ref{primetheorem}) it follows that if $W_a \subset
W_b$ then $A/(b) \to A/(a)$. From proposition \ref{colimite} below it
immediately follows that $E_A$ as constructed in \ref{Spec_A} is
homeomorphic to the output of the Godemment construction applied to
this presheaf (recall Section \ref{terminology}.(\ref{basesobra})).
%(recall Section 1 (\ref{basesobra})).     
\qed\end{sinnadastandard}

\vspace{1ex}

Any prime ideal $P$ with the order relation of $A$ is a filtered poset
(given $a,\, b\, \in P, a \leq a\oplus b,\; b \leq a \oplus b,\; a
\oplus b \in P$), and the assignment $a \mapsto A/(a)$ clearly defines
a directed system indexed by $P$. The following is easy to check:
\begin{proposition}\label{colimite}
Consider any prime ideal $P$, and for each element $a \in P$, the
induced morphism  
$A/(a) \to A/P$. Then, the cone $A/(a) \to A/P$ is a filtered colimit
diagram indexed by $P$. That is, $A/P = \colim{a \in P}{A/(a)}$ 
\qed\end{proposition}

%%%%%%%%%%%%%%%%%%%%%%%%%%%%%%%%%%%%%%%%%%%%%%%%%%%%%%%%%%%%%%%%%%%%%%
\section{Duality between mv-algebras and sheaves of mv-chains} \label{duality}
 %%%%%%%%%%%%%%%%%%%%%%%%%%%%%%%%%%%%%%%%%%%%%%%%%%%%%%%%%%%%%%%%%%%%%%

We introduce now the notion of \emph{mv-Space}. This notion is
inspired in the
notion of \emph{ringed Space} in algebraic geometry (see for example
\cite{H} page 
72). Compare this section with \cite{CDM} section 4.
\begin{definition} {\bf (The category $mv\cc{E}$ of
    mv-spaces)} \label{mv-spaces} 

1) A mv-space is a pair $(X, \, \ff{E})$, where $X$ is a topological
space and \mbox{$\ff{E} = (E \to X)$} is a sheaf of mv-chains on
$X$ (that is, the fibers $E_x,\; x \in X$ are mv-chains). It follows
that for any open set $U \subset X$, the set of sections $\Gamma(U,\, E)$
is a mv-algebra, subalgebra of the product algebra $\prod_{x \in U}
E_x$. We say that $E$ is a sheaf of mv-chains, called the
\emph{structure sheaf}. Technically, $E$ is a mv-chain object in
the topos $\cc{S}h(X)$ of 
sheaves over $X$.

\vspace{1ex}

2) A morphism $(f, \, \varphi): (X, \, \ff{E}) \mr{} (Y, \, \ff{F})$
of mv-spaces is 
a continuous function $X \mr{f} Y$ together with a family $\varphi =
(\varphi_x)_{x \in X}$ of morphisms $F_{f(x)} \mr{\varphi_x} E_x$  such
that for every open $V \subset Y$, and section $s \in \Gamma(V,\, F)$,
the
composite:
$$
\xymatrix
         {
           E_x 
         & F_{f(x)} \ar[l]_{\varphi_x} 
         \\
          x \in  X \ar[u]_{k} \ar[r]^f 
         & Y \ar[u]_s
         }
\hspace{5ex} 
k(x) = (\varphi_x \circ s \circ f)(x).
$$
is a section $k \in \Gamma(U, \, E), \; U = f^{-1}V$. This amounts to
say that $\varphi$ is a morphism of sheaves of mv-chains $f^*F
\mr{\varphi} E$, in the topos $\cc{S}h(X)$, where $f^*$ denotes the inverse
image functor $\cc{S}h(Y) \mr{} \cc{S}h(X)$.

\vspace{1ex}

The reader will have no trouble to check that this defines a category,
the category of mv-spaces, that we denote $mv\cc{E}$.  
\end{definition}

%\begin{sinnadastandard} \label{adjunction1} 
Given any mv-algebra $A$,
  the sheaf of 
  mv-chains $\ff{Spec}_A = (Z_A, \, E_A)$,
\mbox{determines} a mv-space, and this assignment is (contravariantly)
functorial 
$\ff{Spec}: mv\cc{A}^{op} \mr{} mv\cc{E}$. On the other hand, given
any sheaf of \mbox{mv-chains} $\ff{E} = (X,\,E)$, the mv-algebra of
global sections  $\Gamma(\ff{E}) = \Gamma(X, \, E)$, determines a
(contravariant) functor $\Gamma: mv\cc{E}^{op} \mr{} mv\cc{A}$. These
two functors are \emph{\mbox{adjoints} on the right}.
This means that there
is a natural bijection between the hom-sets: 
$$ 
[(X, \, E), \, (Z_A, \, E_A)] \mr{\cong} [A, \,\Gamma(X, \, E)]
$$ 
Under this bijection,
given  $(X, \, E) \mr{(f, \,\varphi)} (Z_A, \, E_A)$, it corresponds
the morphism which to an element $a \in A$, assigns the section $k \in
\Gamma(X, \, E)$ defined by $k(x) = (\varphi_x \circ \widehat{a} \circ
f)(x) = \varphi_x([a]_{f(x)})$.
%\qed\end{sinnadastandard}

By the general theory of adjoint functors (\cite{M}), the
adjunction is equivalent to the
statement in the following proposition, that we leave to check by the
careful reader:

\begin{proposition} \label{adjunction1}
Given any mv-algebra $A$, the following defines a morphism of
\mbox{mv-algebras}, natural in $A$ (see \ref{Spec_A}):
$$
\eta: A
\mr{} \Gamma\ff{Spec}_A = \Gamma(Z_A, E_A), \;\; 
a \mapsto \widehat{a},\;\;  \widehat{a}(P) = [a]_P .
$$
This morphism has the universal property described in the following
diagram: Given any mv-space $(X,\, E)$ and a morphism $A \to \Gamma(X, E)$:
$$
\xymatrix
        {
          A \ar[r]^(.3)\eta \ar[dr]
        & \Gamma(Z_A, E_A) \ar[d]
        \\
        & \Gamma(X, \, E)
        }
\hspace{5ex}
\xymatrix
        {
          (Z_A, E_A) 
        \\
          (X, \, E) \ar@{-->}[u]_{\exists \, !}
        }
$$
\qed\end{proposition}

In the following proposition we describe explicitly the natural
transformation at the other end of the adjunction.

\begin{proposition} \label{adjunction2}
Given any mv-space \mbox{$(X,\, E) \in
mv\cc{E}$,} the unit of the 
adjunction is given by: 
$$
(h,\, \varepsilon) :(X,\, E) \mr{} \ff{Spec}\Gamma(X,\, E) = (Z_{\Gamma(X,\,
  E)}, \,E_{\Gamma(X, \, E)})  
$$
$$
\xymatrix
         {
          E_x \ar@{-}[d]
         & \Gamma(X,\, E)/h(x), \ar[l]_(.65){\varepsilon_x} 
         \\
           x \in X  \ar[r]^h 
         & Z_{\Gamma(X,\, E)} \ar@{-}[u]
         }
\hspace{5ex}
\xymatrix
        { 
         \varepsilon_x([\sigma]_{h(x)}) = \widehat{x}(\sigma) = \sigma(x)
         \\
         h(x) = \{\sigma \, | \, \sigma(x) = 0\}
        }
$$
\end{proposition}
\begin{proof}
The facts concerning this statement are already proved
or straightforward. Remark
that $\varepsilon_x$ is well defined since clearly
$\varepsilon_x(h(x)) = \{0\}$, and $\,h\,$ is continuous since
$h^{-1}(W_\sigma) = [[\sigma = 0]]$ which is open (\ref{sheavesfacts}
(\ref{T_2sheaf})). 
\end{proof}

The general theory of representation of mv-algebras as algebras of
global sections of sheaves of mv-chains that we develop in this paper
is based in the study of the morphism $\eta \,$. 

\vspace{1ex}

\emph{We shall show that for any
mv-algebra, $\,\eta \,$ is an isomorphism.} 

\vspace{1ex}

Again, by the general theory of
adjoint functors, this is equivalent to the fact that the functor
$\ff{Spec}: mv\cc{A}^{op} \mmr{} mv\cc{E}$ is full and faithfull.  As
a consequence, \emph{the dual 
of the category 
of mv-algebras can be considered to be a category of mv-spaces}, which
are objects of a geometrical nature. This is similar to Grothendieck
method in algebraic geometry, which is based in the consideration of the dual
of the category of rings as the category of affine schemes. 

\vspace{1ex}

The fact that $\eta$ is injective is an immediate consequence of 
\cite{COM} 1.2.14 which says that every proper ideal in a mv-algebra is
an intersection of prime ideals. Actually, it is equivalent to the fact that the ideal $\{0\}$ is an intersection of prime ideals.

\begin{proposition}\label{etainjective}
Given any mv-algebra $A$, the morphism  \mbox{$A \mr{\eta} \Gamma(Z_A, E_A)$}
is injective. 
\qed\end{proposition}
\begin{notation} \label{notationeta}
Given any mv-algebra $A$, we shall denote by $\widehat{A}$ the image
of the morphism $\,\eta\,$, thus, $A \mr{\cong} \widehat{A} \subset
\Gamma(Z_A,\, E_A)$
\end{notation}
Given an element $a \in A$, 
the prime ideals of the quotient algebra $A/(a)$ correspond to the prime
ideals in the the set $W_a$ (\cite{COM} 1.2.10). That is,
$Z_{A/(a)} = W_a$. Thus,  the previous proposition applied to the
algebra $A/(a)$ proves the following: 

\begin{proposition} \label{compatible}
Given any mv-algebra $A$ and elements $a,\,b,\,c$ in $A$,
$$
%\hspace{4ex} 
\widehat{b}|_{W_a} = \widehat{c}|_{W_a}  \;\iff\;  [b]_{(a)} =
        [c]_{(a)}\,, \;\;\; that \; is \; \;\; A/(a) \cong \widehat{A}\,|_{W_a}
$$
\qed\end{proposition}

The fact that the ideal $\{0\}$ is an intersection of prime ideals
implies:
\begin{remark} \label{cover}
Given any mv-algebra $A$ and an element $a \in A$, we have $\;Z_A = W_a \;\iff
\; a = 0$. Given any two elements $a, \, b \in A$, $W_a \cup W_b = W_{a 
  \wedge b}$ (see  \ref{mvfacts3} (\ref{p1})). Thus $Z_A = W_a \cup W_b
\; \iff \; a\wedge b = 0$). 
\qed\end{remark}

 In section \ref{posetcoZariski} we shall prove the following lemma:
\begin{lemma}[compacity lemma] \label{compacity}
Given any mv-algebra $A$, the spectral space $Z_A$ is sober, compact,
and has a base of compact open sets.  
\qed\end{lemma}
In particular we know then that $Z_A$ is a compact topological space.   
We analyze now the surjectivity. A global section is determined by a
compatible family of 
sections $\widehat{a_i}$ on an open cover $W_{b_i}$ (see
\ref{sections}).  
We can assume the cover to be finite. Thus, the surjectivity of $\,\eta\,$ amounts to the following:

%\vspace{1ex}
 %
\begin{sinnadaitalica} \label{pushpullp_^n}
Let $a_1, \, \ldots \,a_n  \in A$ be such that $a_1, \, \wedge \, \ldots \, \wedge \, a_n \, = \, 0$. Then, given 
\mbox{$b_1, \, \ldots \,b_n \in A$} such that $(\widehat{b}_i = \widehat{b}_j)|_{Z(\widehat{a}_i) \cap  Z(\widehat{a}_j)}$, there exist a (necessarily unique by \ref{etainjective}) $b \in A$  such that  
$(\widehat{b} = \widehat{b}_i)|_{Z(\widehat{a}_i)}$. 
\end{sinnadaitalica}

 In turn, the validity of \ref{pushpullp_^n} clearly follows by induction from the following: 
\begin{sinnadaitalica} \label{pushpullori}
Let
$a_1,\,a_2 \in A$ be any two elements. Then, given $b_1,\, b_2 \in A$ such
that  $(\widehat{b}_1 = \widehat{b}_2)|_{Z(\widehat{a}_1) \cap  Z(\widehat{a}_2)}$, there
exists $b \in A$,  unique upon restriction to $Z(\widehat{a}_1) \cup  Z(\widehat{a}_2)$, 
such that such that
$(\widehat{b} = \widehat{b_1})|_{Z(\widehat{a}_1}$ and $(\widehat{b} = \widehat{b_2})|_{Z(\widehat{a}_2}$.
This fact means precisely that the following diagram  is a pullback  square:
$$
%\hspace{-2ex} 
%\xymatrix
%       {
%         &  \widehat{A}|_{W_{a_1} \cup W_{a_2} } \ar[r] \ar[d]
%         & \widehat{A}|_{W_{a_1}} \ar[d]
%         %& \hspace{-16ex} \subset \Gamma(W_a,\, E_A) 
%         \\
%         %& \hspace{-13ex} \Gamma(W_b,\, E_A) \supset 
%         &\widehat{A}|_{W_{a_2}} \ar[r]
%         & \widehat{A}|_{W_{a_1} \cap W_{a_2}} 
%           %& \hspace{-6ex} \subset \Gamma(W_a \cap W_b ,\, E_A)
%         }
%\hspace{-3ex}         
\xymatrix
         {
         &  \widehat{A}|_{Z(\widehat{a}_1) \cup Z(\widehat{a}_2)}  \ar[r] \ar[d]
         & \widehat{A}|_{Z(\widehat{a}_1)} \ar[d]
         %& \hspace{-16ex} \subset \Gamma(W_a,\, E_A) 
         \\
         %& \hspace{-13ex} \Gamma(W_b,\, E_A) \supset 
         &\widehat{A}|_{Z(\widehat{a}_2)} \ar[r]
         & \widehat{A}|_{Z(\widehat{a}_1) \cap Z(\widehat{a}_2)} 
         %& \hspace{-6ex} \subset \Gamma(W_a \cap W_b ,\, E_A)
         }
$$
\end{sinnadaitalica}
Since $Z(\widehat{a}_1 \vee  \widehat{a}_2) = Z(\widehat{a}_1) \cap Z(\widehat{a}_2)$ and  $Z(\widehat{a}_1 \wedge \widehat{a}_2) = Z(\widehat{a}_1) \cup Z(\widehat{a}_2)$ it follows that \ref{pushpullori} is equivalent to the following (recall proposition \ref{compatible}):
\begin{lemma}[pushout-pullback lemma] \label{pushpull}
Given any two elements \mbox{$a_1,\,a_2 \in A$,} the following pushout diagram is also a pullback diagram.
$$
\xymatrix
         {
           A/(a_1 \wedge a_2) \ar[r] \ar[d]
         & A/(a_1) \ar[d] 
         \\
           A/(a_2) \ar[r]
         & A/(a_1 \vee a_2)
         }
$$
\qed\end{lemma}
Remark that this diagram is always a pushout since $(a_1 \vee a_2)$ is the supremun
of $(a_1)$ and $(a_2)$ in the lattice of congruences of $A$ (see
\ref{mvfacts1} (\ref{aplusbequalab})). 

In section \ref{pushpulls} we shall prove lemma \ref{pushpull}, thus we have:

\begin{theorem}[Representation Theorem] \label{main}
Given any mv-algebra $A$, the morphism
$$
\eta: A
\mr{} \Gamma\ff{Spec}_A = \Gamma(Z_A, E_A), \;\; 
a \mapsto \widehat{a},\;\;  \widehat{a}(P) = [a]_P .
$$
(see \ref{adjunction1}) is an isomorphism.
\qed\end{theorem}

\section{The maximal and the $[0,\,1]$-valued morphisms
  spectra} \label{msofamvalgebra}

%These algebras are called hypearchimedean although
%there is no notion attached to the name ``archimedean
%mv-algebra''. \emph{We shall simply use ``archimedean'' and not its
%superlative form.} 

\vspace{1ex}

{\sc The maximal spectrum $\ff{Spec}^M_A$}

We consider  now maximal ideals. Recall that maximal ideals are prime.
\begin{sinnadastandard}{ \bf Construction of $\ff{Spec}^M_A$ .} \label{MSpec_A}

The base space $M_A$ is defined to be the subspace $M_A \subset Z_A$
determined by the maximal ideals.  Recall then that a
  base for the topology is given by the sets  $W^M_a = W_a|_{M_A} =
  \{M \; \;|\; M \ni a\}$. The etal space $E^M_A \to M_A$, $E^M_A
  \subset E_A$, is the
  restriction of $E_A$  to $M_A$.
\qed\end{sinnadastandard}

A salient feature of the maximal spectral space is the following:
\begin{proposition} \label{W_a_clopen}
Given any mv-algebra $A$, the sets $W^M_a$ are closed (thus clopen) sets
in $M_A$. 
\end{proposition}
\begin{proof}
We shall show that $\;M_A \backslash W^M_a \;$ is open. Let $P \in M_A$, $P
\notin W^M_a$ . By \ref{mvfacts4} (\ref{maximal}), take an
integer $n$ such that $\neg\, na \in P$. Then $P \in W^M_{\neg\,
  na}$. To finish the proof we have to show that 
$W^M_a \cap  W^M_{\neg\, na} = \emptyset$. 
 We do as follows: Let $Q \in M_A$ be such that $\neg\,
na \in Q$ and $a \in Q$. Then $na \in Q$, thus $\neg\, na \oplus
    na  \in Q$, so $1 \in Q$ (see \ref{mvfacts3} (\ref{p1})), but $Q$ is proper. 
\end{proof}

\begin{corollary} \label{M_AT_2}
Given any mv-algebra $A$, the maximal spectrum  \mbox{$\ff{Spec}^M_A
  \,=\,E^M_A \to M_A$} is a Hausdorff sheaf of simple mv-algebras.   
\end{corollary}
\begin{proof}
Let $P, \, Q \in M_A$, $P \neq Q$. Take $a \in P$, $a \notin Q$. Then
$W^M_a$ and $\;M_A \backslash W^M_a \;$ separate $P$ and $Q$. This
shows that $M_A$ is Hausdorff.

Now let $([a]_P,\,P) \neq ([b]_Q,\,Q)$ in $E^M_A$ (see
\ref{Spec_A}). If \mbox{$P \neq Q$,} take \mbox{$U \ni P$} and
\mbox{$V \ni Q$} be 
disjoint open sets 
in $M_A$. Then  \mbox{$\widehat{a}(U) \ni ([a]_P,\,P)$} and
\mbox{$\widehat{b}(V)\ni ([b]_Q,\,Q)$} are disjoint open sets in
$E_A$. If $P = Q$ \mbox{(denote both $R$),} then \mbox{$[a]_R \neq [b]_R$.}
Consider the closed 
set in $M_A$,  
$W^M_{d(a,\,b)}= [[\widehat{a} = \widehat{b}]]$ (see \ref{notation}). Notice
%$ \mbox{$\{R \,| \; [a]_R = [b]_R \} =$} \mbox{$ \{R \, |
%\; \widehat{a}(R) = \widehat{b}(R)\}$}
that \mbox{$R \notin  W^M_{d(a,\,b)}$.}  Then $\widehat{a}(M_A
\backslash W^M_{d(a,\,b)}) \ni ([a]_R,\,R) $ and
\mbox{$\,\widehat{b}(M_A \backslash 
W^M_{d(a,\,b)}) \ni ([b]_R,\,R)$} are disjoint open sets in
$E^M_A$. This finishes the proof that  $E^M_A$ is Hausdorff.

Finally, by \ref{mvfacts4} (\ref{maximal2}), we know that the fibers
are simple mv-algebras.   
\end{proof} 

Given a prime ideal $P$, there exist a
maximal ideal $M \supset P$, \ref{mvfacts4} (\ref{maximal}). Thus, given any 
open set $W_a$ in the base of the topology, we have that $P \in W_a \;
\rimply  M \in W_a$. It follows:
\begin{comment} \label{onlyprime}
The spectral space $Z_A$ is never Hausdorff unless it is equal to $M_A$.
\qed\end{comment}
\begin{comment} \label{dense}
The inclusion $M_A \subset Z_A$ is dense. Namely, the prime spectral space
$Z_A$ is the closure of the maximal spectral space $M_A$. 
\qed\end{comment}

We study now the injectivity of the morphism $A \mr{\eta} \Gamma(M_A, E^M_A)$, $a \mapsto  \widehat{a}$, $\widehat{a}(M) = [a]_M$.
Recall that a mv-algebra $A$ is \emph{semisimple} if it is non trivial and the
intersection of all maximal ideals is $\{0\}$. Thus, like proposition \ref{etainjective} we have:
%
 %equivalently, if the
%implication  \mbox{($\chi(a) = 
%\chi(b) \;\; \forall 
%\chi \in X_A \;\;\rimply\;\; a = b$)} holds. 
%
\begin{proposition}\label{etaMinjective}
A mv-algebra is semisimple if and only the morphism $A \mr{\eta} \Gamma(M_A, E^M_A)$ is injective. 
\qed\end{proposition}
Given an element $a \in A$, 
the maximal ideals of the quotient algebra $A/(a)$ correspond to the maximal
ideals in the the set $W^M_a$ (\ref{mvfacts4} (\ref{phimenosunomaximal})). That is,
$M_{A/(a)} = W^M_a$. Now we need a proposition for the maximal spectrum corresponding to proposition \ref{compatible}. This is not possible for semisimple algebras in general, a stronger hypothesis on the algebra is necessary.
\begin{definition} \label{ssemisimple}
A mv-algebra $A$ is \emph{strongly semisimple} if $A$ and all its quotients $A/(a)$, $a \in A$, are semisimple.
\end{definition}
Clearly if $A$ is strongly semisimple, so are all its quotients $A/(a)$. \mbox{Hyperarchimedean} algebras (section \ref{hyperarchimedean}) are strongly semisimple. As we shall see, finitely presented algebras (section \ref{fp}) are examples of non hyperarchimedean strongly semisimple algebras. The example following corollary 3.4.4 in \cite{COM} shows that the semisimple mv-algebra $Cont([0,\,1],\, [0,\,1])$ is not strongly semisimple. The following is immediate:

\begin{proposition} \label{compatibleM}
A mv-algebra $A$ is strongly semisimple if and only if given elements $a,\,b,\,c$ in $A$,
$$
\widehat{b}|_{W^M_a} = \widehat{c}|_{W^M_a}  \;\iff\;  [b]_{(a)} =
        [c]_{(a)}\,, \;\;\; that \; is \; \;\; A/(a) \cong \widehat{A}\,|_{W^M_a}
$$
\qed\end{proposition}
The reader will have no difficulty proving the following:
\begin{corollary} \label{strongsemi}
Let $A$ be a strongly semisimple mv-algebra, then the restriction morphism
$
\Gamma(Z_A, E_A) \mr{\rho}  \Gamma(M_A, E^M_A)
$
is injective. 
\qed\end{corollary}

Notice that from the representation theorem \ref{main} it immediately follows that the statement in this corollary holds for any semisimple algebra. However, it is worth to pay the price of the stronger hypothesis to have this simple proof independent of the representation theorem. An interesting application of this corollary is proposition \ref{strongsemi2} and theorem \ref{mc=gamma}.

\vspace{1ex}

{\sc The [0,\,1]-valued morphisms spectral space}

\vspace{1ex}

We consider
now a different construction, that we denote $X_A$, of a topological space
associated to a mv-algebra $A$. This construction is more akin to
functional analysis than to algebraic geometry, and has been
considered in particular in \cite{CDM}.

\vspace{1ex}

\begin{sinnadastandard}{ \bf Construction of $X_A$ .} \label{X_A} 
 The points of $X_A$ are all the $[0,\,1]$-valued morphisms,
 \mbox{$X_A = [A,\, [0,\,1]\,] = \{\chi: A  \mr{} [0,\, 1]\} \subset [0,
     \, 1]^A $.}  The 
 topology is  
  the subspace topology inherited from the product space, with the unit
  interval endowed with the usual topology. Recall then that a subbase
  for this topology is given by the sets $W^X_{a,\,U} = \{\chi \, | \;
  \chi(a) \in U \}$, for $a \in A$ and $U$ an open set,
  $U \subset [0,\, 1]$. 
\qed\end{sinnadastandard} 
 
\begin{remark} \label{WXaclosed}
Given any $x \in [0,\,1]$,  
 \mbox{$W^X_{a,\,x} = \{\chi \, | \;  \chi(a)= x  \}$} is a closed set
of $X_A$. When $x = 0$ we write $W^X_a = \{\chi \, | \;  \chi(a)= 0  \}$.
\qed\end{remark}

\emph{Since $X_A$ is closed in $[0, \, 1]^A$, $X_A$ is
  a compact Hausdorff space}.

\vspace{1ex}

From \ref{mvfacts4} (\ref{uniquesubalgebra}) it follows:

\begin{proposition} \label{X=M}
Maximal ideals $M$ are in bijection with morphisms \mbox{$A
\mr{\chi} [0,\, 1]$.} % = \{a\,|\, \chi(a) = 0\}$.
If $M$ corresponds to $\chi$, $A/M \mr{\cong} \chi(A) \subset
[0,\,1]$, and \mbox{$\; a \in M \;\iff\;
  \chi(a) = 0$.}
\qed\end{proposition}

We denote $X_A \mr{\kappa} Z_A$ the injection defined by
\mbox{$\kappa(\chi) = ker(\chi)$.} Its image is the subspace $M_A
\subset Z_A$ of maximal ideals. We shall (abuse notation) and denote $\,\chi\,$
the inverse map $\chi = \kappa^{-1} : M_A \mr{} X_A$, and write
$\chi(M) = \chi_M$. 

\vspace{1ex}

Each
element $a \in A$ defines a continuous function  
$X_A \mr{\widehat{a}} [0,\,1]$, \mbox{$\widehat{a}(\chi) = \chi(a)$.}
Thus,  $W^X_{a,\,U} = \widehat{a}^{\,-1}U$. 
By definition, \emph{the topology in $X_A$ is
the initial topology with respect to the functions $\widehat{a}$, $a \in A$.}

\vspace{1ex}

For any topological space $X$ we denote by $Cont(X,\,[0,\,1])$, the
\mbox{mv-algebra} of all $[0,\,1]$-valued continuous functions on X.
 Recall that a subalgebra \mbox{$A \subset
Con(X, \, [0,\,1])$} is said to be \emph{separating} if given any two
points $\chi \neq
\xi$ in $X$, there is  $f \in A$ such that $f(\chi) = 0$ 
and $f(\xi) > 0$.

\vspace{1ex}

If two $[0,\,1]$-valued
morphisms are different, by proposition \ref{X=M} their kernels must
be different.  It follows:
\begin{proposition} \label{separating}
Given any mv-algebra $A$, the
functions $X_A \to [0,\,1]$ of the form $\,\widehat{a}\,$, $a \in A$, form a
separating subalgebra of the mv-algebra
$Cont(X_A ,\; [0,\,1])$. That is, given any two points $\chi \neq
\xi$ in $X_A$, there is  $a \in A$ such that $\chi(a) = 0$ 
and $\xi(a) > 0$.
\qed\end{proposition}

\vspace{1ex}

Given $a \in A$, we denote $W^X_a \,=\, \{\chi \, | \;  \chi(a)= 0  \} \,=\, 
X_A \backslash  W^X_{a,\, (0,\,1]}$. Notice then that
\mbox{$\kappa(W^X_a) = W^M_a$} and $\chi(W^M_a) = W^X_a$.

\begin{proposition} \label{M_AtoX_A}
Given any mv-algebra $A$, the sets $W^X_{a,\, (0,\,1]} = $
  \mbox{$ = \{\chi \;|\; 
  \widehat{a}(\chi)  > 0\} =  \{\chi \;|\; 
  \widehat{a}(\chi)  \neq 0\}$} form a base for the topology of $X_A$.
\end{proposition}
\begin{proof}
 Since $X_A$ is a compact Hausdorff
space, the result follows from  proposition \ref{separating} and
(\cite{COM}, Remark to Theorem 3.4.3).
\end{proof}

\emph{Thus, the open sets basis $W^X_{a,\, (0,\,1]}$  is formed by the
complements of the Zero sets of the functions $\widehat{a}$.}

\vspace{1ex}

Under the bijection $M_A \cong X_A$ determined by $(\kappa, \, \chi)$,
we have $\kappa(W^X_{a,\, (0,\,1]}) = \{M \in M_A \;|\; a \notin
  M\}$. So the topology of $X_A$ corresponds to the
  \mbox{\emph{Zariski} topology} in $M_A$, while 
  the topology of $M_A$ was defined to be the \mbox{\emph{coZariski}
    topology.}  

From proposition \ref{W_a_clopen} it follows:

\begin{proposition} \label{M_AfinerthanX_A}
Given any mv-algebra $A$, the bijection $M_A \mr{\chi} X_A$ is
continuous. That is, in the set $M_A$ of 
maximal ideals the coZariski topology is finer than the Zariski topology.
\qed\end{proposition}

\begin{remark}\label{M_Anotcompact}
The maximal spectral space $M_A$ in not a compact space in general.
It is compact if and only if it is homeomorphic to $X_A$ under the
bijection $\, \chi \,$.
\end{remark}
\begin{proof}
One implication is clear. For the other, assume $M_A$ to be
compact. Then the continuous map $M_A \mr{\chi} X_A$ 
is also a closed map (thus an homeomorphism) because $X_A$ is Hausdorff. 
\end{proof}

In general the coZarisky topology will be strictly finer than the
Zariski topology. In fact, the
topologies coincide for 
\emph{hyperarchimedean} mv-algebras, and only in this case (see
\mbox{section \ref{hyperarchimedean}).}
%We finish thissection with: 

\begin{proposition} \label{M_AhomeoX_A}
The map (injection) $X_A \mr{\kappa} Z_A$ is continuous
if and only if the sets $W^X_a = \{\chi \,|\; \chi(a) = 
0\}$ are open (thus clopen) sets in $X_A$. If such is the case, then
$\, \kappa \,$ and $\, \chi \,$
  establish an homeomorphism $X_A \mr{\cong} M_A$.
\end{proposition}
\begin{proof} 
$W^X_a = \kappa^{-1}W_A$, which are an open basis of $Z_A$, thus
  $\kappa$  
  is continuous if and only if the sets $W^X_a$ are open in
  $X_A$. In this case then, the continuous bijection $\,\chi\,$ (proposition
  \ref{M_AfinerthanX_A}) has a continuous inverse.
  The proof also
  can be completed without using this proposition. In fact, if
  $\kappa$ is continuous, 
  it is a closed map as a map $X_A \mr{\kappa} M_A$, because $X_A$ is 
  compact and $M_A$ is Hausdorff. 
\end{proof}

The mv-algebra of global sections of  $\ff{Spec}^M_A 
\,=\,E^M_A \to M_A$ is  related with the mv-algebra of continuous functions
$Cont(M_A, \, [0,\,1])$. Define a map $E^M_A \mr{\lambda} [0,\,1]$
by writing $\lambda([a]_M ,\, M) = \chi_M(a))$ for each $([a]_M ,\, M)
\in E_M$. Given any $a\in A$, observe that the following diagram is
commutative:
$$
\xymatrix
         {
            E^M_A \ar[r]^{\lambda} 
          & [0,\,1]
          \\
            M_A \ar[u]^{\widehat{a}} \ar[r]^{\chi}
          & X_A \ar[u]^{\widehat{a}}
         }
$$
Since the topology of $E^M_A$ is the final topology with respect to
the functions $\widehat{a}(M) = ([a]_M,\,M)$, and the topology of $X_A$ is
the initial topology with respect to the functions $\widehat{a}(\chi)
= \chi(a)$, we have that  \mbox{\emph{$\lambda$ is continuous
    $\,\iff\,$ $\chi$ is 
  continuous}.} It follows then from proposition \ref{M_AfinerthanX_A}
that the function $\lambda$ is a continuous function. The whole
statement concerning this fact is the following:
\begin{proposition} \label{GammasubCon}
Given any mv-algebra $A$, the map $E^M_A \mr{\lambda} [0,\,1]$ defined
by $\lambda([a]_M ,\, M) 
= \chi_M(a)$ is continuous and establishes:

%\vspace{10ex}

a) A continuous injection 
$$E^M_A \mmr{} [0,\,1] \times M_A \;\; over \;\; M_A\,,$$
where $[0,\,1]$ is endowed with the usual topology. 

b) By composition an injective morphism 
$$\Gamma\ff{Spec}^M_A = \Gamma(M_A, E^M_A) \mmr{\lambda_*}
Cont(M_A,\,[0,\,1])\,, \;\;
\lambda_*(\widehat{a})(M) = \chi_M(a).
$$
Furthermore, this morphism establishes an isomorphism of
$\Gamma\ff{Spec}^M_A$ onto a 
separating subalgebra of $Cont(M_A,\,[0,\,1])$.
\end{proposition}
\begin{proof}
Besides \ref{M_AfinerthanX_A}, proposition
\ref{separating} completes the proof. It is also a particular case of
proposition \ref{Gammaisarchi} below. 
\end{proof}
Consider now a global section $M_A \mr{\sigma} E^M_A$ which is the restriction of a global section of the prime spectrum $Z_A \supset  M_A$. Since $Z_A$ is a compact topological space (corollary \ref{Z_Acompact}), there is a finite open cover $W^M_{a_i} \subset M_A$ over which $\sigma$ is of the form $\widehat{b}_i$ for some $b_i \in A$. Thus the composite \mbox{$X_A \mr{\kappa} M_A \mr{\sigma} E^M_A$} is of the form $\widehat{b}_i$ over a finite cover of $X_A$ by the closed sets $W^X_{a_i}$ (remark \ref{WXaclosed}). It follows then that the composite $X_A \mr{\sigma \kappa} E^M_A$ is  a continuous function, thus so is the composite $X_A \mr{\lambda \sigma \kappa} [0,\,1]$. This determines a morphism 
\begin{equation} \label{kappalambda}
\Gamma(Z_A, E_A) \mr{\kappa^*\lambda_*}  Cont(X_A, \, [0,\,1])
\end{equation} 
$$\;\;\kappa^*\lambda_*(\sigma) = \lambda \sigma \kappa\,, \;\; 
\lambda \,\widehat{a} \,\kappa(\chi) = \chi(a).
$$
which makes the following diagram commutative (recall \mbox{$\kappa\chi = id$}):
$$
\xymatrix
             {
               \Gamma(Z_A, E_A) \ar[r]^{\kappa^*\lambda_*} 
                                             \ar[d]^{\rho}
             &
               Cont(X_A, \, [0,\,1]) \ar@{^{(}->}[d]^{\chi^*}
             \\ 
               \Gamma(M_A,\, E^M_A) \ar@{^{(}->}[r]^{\lambda_*} 
              &
               Cont(M_A,\,[0,\,1])
             }                             
$$
%
%\ell(\sigma)(\chi) = \lambda \sigma \kappa(\chi) = \lambda_*(\sigma)
%(\kappa(\chi)) = \chi(a)$ (where $\sigma(P) = ([a]_P,\,P)$).
%
It follows (recall corollary \ref{strongsemi}):
\begin{proposition} \label{strongsemi2}
Let $A$ be a strongly semisimple mv-algebra (\ref{ssemisimple}). Then the morphism in \ref{kappalambda} is injective. 
$
\xymatrix
             {
               \Gamma(Z_A, E_A) \ar@{^{(}->}[r]^{\kappa^*\lambda_*} 
            & 
               Cont(X_A, \, [0,\,1])
            }
$
\qed\end{proposition}

\vspace{1ex}

%%%%%%%%%%%%%%%%%%%%%%%%%%%%%%%%%%%%%%%%%%%%%%%%%%%%%%%%%%%%%%%%%%%%%%
\section{Hyperachimedian algebras} \label{hyperarchimedean}
 %%%%%%%%%%%%%%%%%%%%%%%%%%%%%%%%%%%%%%%%%%%%%%%%%%%%%%%%%%%%%%%%%%%%%%

Before attacking the representation theorem in the general case we
consider the particular case of hyperarchimedean algebras
\cite{CDM}, \cite{COM}. 

%These algebras are called hypearchimedean although
%there is no notion attached to the name ``archimedean
%mv-algebra''. \emph{We shall simply use ``archimedean'' and not its
%superlative form.} 

For hyperarchimedean mv-algebras every prime ideal is
maximal, the Zariski and coZariski topologies coincide, and the spectrum space
$Z_A$ is separated (Hausdorff).
Furthermore, all the mv-chains $A/P$ are subalgebras of the 
real unit interval $[0,\, 1]$. All these facts allow us to apply classical
methods as in \cite{CDM}. 
%
%In fact, the mv-spaces which come into consideration have as structure
%sheaf the sheaf of germs of continuous $[0,\, 1]$-valued functions.
%
\vspace{1ex}

We shall generalize the representation of locally finite
mv-algebras  by Cignoli-Dubuc-Mundici in
\cite{CDM} to general hyperarchimedean mv-algebras.

\vspace{1ex}

The notion of hyperarchimedean mv-algebra is intimately related with the
notion of
maximal ideal. 
Hyperarchimedean algebras are semisimple algebras such that all its
quotients are also semisimple (\cite{COM} 6.3.2). Recall:

\begin{definition}[\cite{COM} 6.2.4]  \label{archidef} $ $

a) An element $a \in A$ in a mv-algebra is \emph{archimedean} if there is
an integer $n \geq 1$ such that $na = (n+1)a$.

b) A mv-algebra A is \emph{hyperarchimedean} if every $a \in A$ is 
archimedean.
\end{definition}

\begin{remark}\label{archi1} 
Clearly the real unit interval is an
hyperarchimedean \mbox{mv-algebra} (\ref{mvfacts1} (\ref{nx})(\ref{[01]})), and if $x \in [0, \, 1]$ is such that $nx = (n+1)x$, then
$x \geq 1/n$. 
\qed\end{remark}
\begin{proposition}[\cite{COM} 6.3.2] \label{archiprop} $ $

a) A mv-algebra is hyperarchimedean if and
  only if every prime ideal is maximal. 

b) A mv-chain is hyperarchimedean if
  and only if is a subalgebra of the real unit interval $[0,\, 1]$.
\qed\end{proposition}

We take from \cite{CDM} 4.5 the argument in the proof of the following
proposition: 

\begin{proposition} \label{X_AStone}
Given an archimedean element $a \in A$ in a mv-algebra $A$, the set 
$W^X_a = \{\chi \, | \; \chi(a) = 0\}$ is open in $X_A$, thus a clopen
set.
\end{proposition}
\begin{proof}
Suppose that $a \in A$ is hyperarchimedean. Take an integer $n \geq 1$ such
that $na = (n+1)a$, whence $X_A \backslash W^X_a = \{\chi \, | \;
\chi(a) > 0\} =  
\{\chi \, | \; \chi(a) \geq 1/n \}$ (remark \ref{archi1}), which is
closed, thus $W^X_a$ is open. 
\end{proof}

\begin{observation} \label{reflectsarchi}
Given a semisimple mv-algebra $A$, and  $a \in A$, if
$\chi(a)$ is archimedean with \emph{the same} n for all $\chi \in
X_A$, then $a$ is archimedean. Actually, it is enough that 
the assumption holds for all $\chi \notin  W^X_a$.
\qed\end{observation}
For semisimple mv-algebras the converse of proposition \ref{X_AStone}
is also valid.
\begin{proposition} \label{alsoconverse}
Given a semisimple mv-algebra $A$, if the set $W^X_a$ is open in
$X_A$, then the element $a$ is archimedean.
\end{proposition}
\begin{proof}
Suppose that $W^X_a$ is open. Then $X_A \backslash W^X_a$
is closed, thus compact. It follows that the set $\{\chi(a)\,|\; \chi
\in X_A \backslash W^X_a\}$ is separated from $0$. Let $n \geq 1$ be such that
$X_A \backslash W^X_a = \{\chi \, | \; \chi(a) \geq 1/n \}$. Then for
all $\chi$ in $X_A \backslash W^X_a$, $n\chi(a) = (n+1)\chi(a)$. The
proof finishes by \ref{reflectsarchi}.
\end{proof}
We list now a series of conditions that characterize hyperarchimedean
\mbox{mv-algebras.}

\begin{theorem} \label{archiequivalences} $ $
%\begin{enumerate}

a) The following conditions in an arbitrary mv-algebra $A$ are equivalent:
\begin{enumerate}
\item $A$ is hyperarchimedean. \label{1}
\item $Z_A = M_A$ (that is, every prime ideal is maximal). \label{2}
\item The prime spectral space $Z_A$ is Hausdorff. \label{4}
\item The prime spectrum $\ff{Spec}_A = (E_A \to Z_A)$ is a Hausdorff
 sheaf. \label{5} 
\item For all $a \in A$, the sets $W_a \subset Z_A$ are closed (thus
  clopen) in $Z_A$. \label{6}
\item  The map $\kappa: X_A \mr{\cong} Z_A$ is an
  homeomorphism (in particular $X_A$ is
  homeomorphic to $M_A$, $\kappa: X_A 
  \mr{\cong} M_A$). \label{9}
\item The prime spectral space $Z_A$ is a Stone space (compact Hausdorff
  with an open base of clopen sets). \label{10}

\vspace{1ex}

\hspace{-9ex}  b) The following conditions in a semisimple mv-algebra
$A$ are equivalent:

\vspace{1ex} 

\hspace{-7.5ex} $(1)$ $A$ is hyperarchimedean.

\item For all $a \in A$, the sets $W^X_a \subset
  X_A$ are open (thus clopen) in $X_A$. \label{7}
\item The map $\kappa: X_A \mr{} Z_A$ is continuous. \label{8} 
\item The maximal spectral space $M_A$ is compact. \label{3}
\item The map  $\kappa: X_A \mr{} M_A$ is an homeomorphism. \label{11}
\end{enumerate}
\end{theorem}
\begin{proof} $ $

a)

(\ref{1}) $\iff$ (\ref{2}): By \ref{archiprop}. 

(\ref{2}) $\iff$ (\ref{4}) and (\ref{2}) $\iff$ (\ref{5}): By \ref{M_AT_2}, \ref{onlyprime}.

(\ref{2}) $\rimply$ (\ref{6}): By \ref{W_a_clopen}.

(\ref{6}) $\rimply$ (\ref{4}): Let $P,\, Q \in Z_A$, $P \neq Q$. Take $a \in P$,
  $a \notin Q$. Then, $W_a$ and $Z_A\backslash W_a$ are open sets that
  separate $P$ and $Q$.

(\ref{1}) $\rimply$ (\ref{9}): By \ref{M_AhomeoX_A} and \ref{X_AStone}.

(\ref{9}) $\rimply$ (\ref{10}): Recall that $X_A$ is a
  compact Hausdorff space. 

(\ref{10}) $\rimply$ (\ref{4}): A Stone space, in particular, is
  Hausdorff.

\vspace{1ex}

b)

(\ref{1}) $\rimply$ (\ref{7}): By \ref{X_AStone}.

(\ref{7}) $\rimply$ (\ref{1}): By \ref{alsoconverse}.

(\ref{7}) $\iff$ (\ref{8}): Notice that $W^X_a = \kappa^{-1}W_a$.

(\ref{8}) $\rimply$ (\ref{11}): By \ref{M_AhomeoX_A}.

(\ref{11}) $\rimply$ (\ref{8}): Since the inclusion $M_A \subset Z_A$
is continuous.

(\ref{11}) $\iff$ (\ref{3}): By \ref{M_Anotcompact}
\end{proof}

 From theorem
\ref{archiequivalences} it
clearly follows:

\begin{theorem} \label{archistone}
A mv-algebra $A$ is hyperarchimedean if and only if the spectrum sheaf
$\,\ff{Spec}_A$ is a Hausdorff sheaf of simple mv-algebras over a Stone space.
\qed\end{theorem}

Let $\ff{E} = (E \to X)$ be any  Hausdorff sheaf of simple
mv-algebras. Consider the mv-algebra of global
sections $\Gamma\ff{E} = \Gamma(X, \, E)$. 
The following proposition has been proved, although not explicitly,
in \cite{CDM} section 6.  
\begin{proposition} \label{Gammaisarchi}
For each $\sigma \in
\Gamma(X, E)$ define $f_\sigma(x) = \lambda_x(\sigma(x))$, where
$\lambda_x$ is the unique embedding $E_x \mmr{} [0, \,1]$
\mbox{(\ref{mvfacts4} (\ref{uniquesubalgebra})).} The assignment
$\sigma \mapsto f_\sigma$ defines an 
embedding $\Gamma(X, E) \mmr{} Con(X, \, [0,\,1])$ into an
hyperarchimedean subalgebra of $Con(X, \, [0,\,1])$. It is separating
if $X$ has a base of clopen sets. 
\end{proposition}
\begin{proof}
In
\cite{CDM} 6.4 it is proved that the function $f_\sigma$ is
continuous (here it is used that $E$ is Hausdorff). Since the zero set
\mbox{$Z(f_\sigma) = [[\sigma = 0]] \subset X$} is open
\mbox{(\ref{sheavesfacts} (\ref{T_2sheaf})),} it follows, by 
\cite{CDM} 4.5, that each $f_\sigma$ is an archimedean
element. Finally the last assertion follows because the characteristic
function (as a section) of any clopen set is a continuous section.
\end{proof}

%We say that a sheaf  $\ff{E} = (E \to X)$ is \emph{global} if each $e
%\in E$ is the image of some global section. 

Global Hausdorff sheaves of simple mv-algebras over Stone spaces are
completely determined by their algebra of global sections.
\begin{theorem} \label{epsiloniso}
Given any  mv-space $(X, \ff{E})$, with $\ff{E}$ a global sheaf, $X$ a
Stone space and $E$
Hausdorff, then the unit of the adjunction (see \ref{adjunction2}):
$$
(h,\, \varepsilon):(X,\, E) \mr{} \ff{Spec}\Gamma(X,\, E) = (Z_{\Gamma(X,\,
  E)}, \,E_{\Gamma(X, \, E)})  
$$
is an isomorphism of mv-spaces. 
\end{theorem}
\begin{proof}

We refer to proposition \ref{adjunction2}. Consider the map 
\mbox{$X \mr{h} Z_{\Gamma(X,\, E)}$,} \mbox{$h(x) = \{\sigma \, | \, \sigma(x) =
0\}$.} By the previous proposition we know that \mbox{$Z_{\Gamma(X,\, E)} =
M_{\Gamma(X,\, E)}$.} The ideal $h(x)$ corresponds to the ideal
$\{f_\sigma \, | \, f_\sigma (x)= 0\}$. Then, \cite{COM} theorem 3.4.3
(i) shows that $\, h \,$ is a bijection. Since it is continuous, it
follows that it is an homeomorphism because $X$ is compact and
$M_{\Gamma(X,\, E)}$ is Hausdorff. Finally, the morphism
$\,\varepsilon_x\,$ is injective since the algebras are 
simple, and surjective by definition of global sheaf.  
\end{proof}
From proposition \ref{Gammaisarchi} and theorem \ref{epsiloniso} it
immediately follows a companion of theorem \ref{archistone}.
\begin{theorem} \label{stonearchi}
A sheaf $\ff{E} = (E \to X)$ is a global Hausdorff sheaf of simple
mv-algebras over a Stone space if and only if the mv-algebra of global
sections $\Gamma\ff{E} = \Gamma(X, \, E)$ is hyperarchimedean.
\qed\end{theorem}

\vspace{1ex}

%\begin{sinnadastandard} \label{pushpullsection}
{\sc The pushout-pullback lemma  for hyperarchimedean algebras.}
 %\end{sinnadastandard} 

\vspace{1ex}

We prove now lemma \ref{pushpull} for hyperarchimedean algebras.
 Recall that in this case $Z_A = M_A$, and that by proposition \ref{GammasubCon} b) we have an injective morphism 
$\widehat{A}  \mmr{\lambda_*} Cont(M_A,\,[0,\,1])$.

The pushout-pullback lemma can be stated as follows (see  \ref{pushpullori}):

\begin{lemma}  \label{pushpullhyperarchimedean}
Given any hyperarchimedean algebra $A$, and any two elements $a_1,\,a_2 \in A$, the following diagram of $[0,\,1]$ valued functions is a pullback. 
$$
\xymatrix
         {
         &  \lambda_*(\widehat{A})|_{W_{a_1} \cup W_{a_2} } \ar[r] \ar[d]
         & \lambda_*(\widehat{A})|_{W_{a_1}} \ar[d]
         %& \hspace{-16ex} \subset \Gamma(W_a,\, E_A) 
         \\
         %& \hspace{-13ex} \Gamma(W_b,\, E_A) \supset 
         &\lambda_*(\widehat{A})|_{W_{a_2}} \ar[r]
         & \lambda_*(\widehat{A})|_{W_{a_1} \cap W_{a_2}}
         %& \hspace{-6ex} \subset \Gamma(W_a \cap W_b ,\, E_A)
         }
$$
\end{lemma}
\begin{proof}
 Let $b_1,\, b_2 \in A$ such
that $(\lambda_*\widehat{b_1})|_{W_{a_1}}$,
$(\lambda_*\widehat{b_2})|_{W_{a_2}}$ are compatible in the intersection
\mbox{$W_{a_1} \cap W_{a_2}$.} We have to show that there 
exists \mbox{$b \in A$}, unique upon restriction of $\lambda_*\widehat{b}$ to $W_{a_1} \cup W_{a_2}$, 
such that such that
\mbox{$(\lambda_*\widehat{b})|_{W_{a_1}} = (\lambda_*\widehat{b}_1)|_{W_{a_1}}$} and $(\lambda_*\widehat{b})|_{W_{a_2}} =
(\lambda_*\widehat{b}_2)|_{W_{a_2}}$. 

\vspace{1ex}

Take an integer $n \geq 1$ such that
$na_1 = (n+1)a_1$, $na_2 = 
(n+1)a_2$. Then 
$$
b = (b_1 \wedge b_2) \vee (na_1 \wedge b_2) \vee (na_2 \wedge b_1)
$$
is the required element. 

To check this first we simplify notation: Let $f_1 = \lambda_*\widehat{a_1}$, 
$f_2 = \lambda_*\widehat{a_2}$, $g_1 = \lambda_*\widehat{b_1}$, $g_2 = \lambda_*\widehat{b_2}$, and $g = \lambda_*\widehat{b}$.
Since $\eta$ and $\lambda_*$ are morphisms, we have $nf_1 = (n+1)f_1$, $nf_2 = (n+1)f_2$ and  
$$
g = (g_1 \wedge g_2) \vee (nf_1 \wedge g_2) \vee (nf_2 \wedge g_1)
$$
Remark also that $W_{a_1} = Z(f_1)$ and $W_{a_2} = Z(f_2)$. 

%Given  any $x \in M_A$, 
%$
%f_1(x) > 0 \; \rimply   nf_1(x) = 1, \;\; 
%f_2(x) > 0 \; \rimply   nf_2(x) = 1.
%$

%
%$$
%lambda_*\widehat{b} = (\lambda_*\widehat{b_1} \wedge \lambda_*\w%idehat{b_2}) \vee (n\,\lambda_*\widehat{a_1} \wedge \lambda_*\
%$widehat{b_2}) \vee (n\,\lambda_*\widehat{a_2} \wedge \lambda_*
%\widehat{b_1})
%$$
%

\vspace{1ex}

%$$
%\lambda_*\widehat{a_i})(M) > 0 \;\;\;\; \rimply  \;\;\; (n\,\lambda_*
%\widehat{a_i})(M)  = 1, \hspace{4ex}   i =1,\,2 .
%$$ 
Let $x \in M_A$. Then $x \in Z(f_1)$ or $x \in Z(f_2)$. Suppose  $x \in
 Z(f_1)$,  thus $ f_1(x) = 0$. Consider two cases: 

\vspace{1ex}

$ x \notin  Z(f_2), \; f_2(x) > 0$, thus  $nf_2(x) = 1$. Then:

\hfill $ g(x)  = (g_1(x) \wedge
g_2(x)) \vee (g_1(x) = g_1(x).
$

\vspace{1ex}

$
x \in Z(f_2), \; f_2(x) = 0$ and  $g_1(x) = g_2(x)$. Then:

\hfill $g(x) = g_1(x) \wedge g_2(x) = g_1(x)$.

\vspace{1ex}

For $x \in Z(f_2)$ we do in
the same way. 
\end{proof} 

We adapt now a formula taken from the proof of  \cite{NPM} lemma 5.3 to remark the following:
\begin{remark} \label{putaputa}
In the previous proposition, the formula 
$$
c = (b_1 \odot  \neg na_1) \vee (b_2 \odot \neg na_2)
$$
also gives the required element.
\end{remark} 
\begin{proof}
We let the reader to check this in the same way that for the formula $b$ in the proof of proposition \ref{pushpullhyperarchimedean}.
\end{proof}
 By theorem \ref{archiequivalences} (\ref{10}) we know, in particular,
 that lemma \ref{compacity} holds for hyperarchimedean algebras. Thus
 we have:
\begin{theorem} \label{etaiso}
Given any  hyperarchimedean mv-algebra $A$, the unit of the adjunction
(see \ref{adjunction1}): 

$
\hspace{18ex} \eta: A \mr{\eta} \Gamma\ff{Spec}_A = \Gamma(Z_A, E_A).
$

\noindent is an isomorphism of mv-algebras.
\qed\end{theorem}

Algebras of global sections of global Hausdorff sheaves over Stone
spaces are known as \emph{boolean products}. Remark then that this
theorem together with theorem \ref{archiequivalences} 
(\ref{10}) yields the characterization of hyperarchimedean mv-algebras
as boolean products of simple mv-algebras, \cite{COM} 6.5.6.

\vspace{1ex}

Putting together theorems \ref{epsiloniso} and \ref{etaiso} we have
(see section \ref{duality}):
\begin{theorem} \label{equivalence}
The functors $\; mv\cc{A}^{op} \mr{\ff{Spec}} mv\cc{E}\;$ and 
$\; mv\cc{E}^{op} \mr{\Gamma} mv\cc{A}$ establish an equivalence
between the dual 
of the category of hyperarchimedean mv-algebras and the category of
mv-spaces which are global Hausdorff sheaves over Stone spaces.  
\qed\end{theorem}

In view of theorem \ref{archiequivalences} (\ref{9}), this
theorem generalizes (and yields) the equivalence for the category of
locally finite mv-algebras established in \cite{CDM} \mbox{theorem 6.9.} 

%%%%%%%%%%%%%%%%%%%%%%%%%%%%%%%%%%%%%%%%%%%%%%%%%%%%%%%%%%%%%%%%%%%%%%%
\section{Finitely presented mv-algebras} \label{fp}
%%%%%%%%%%%%%%%%%%%%%%%%%%%%%%%%%%%%%%%%%%%%%%%%%%%%%%%%%%%%%%%%%%%%%%%

 In this section we consider finitely presented mv-algebras and prove the pushout-pullback lemma in this case. A
mv-algebra is of \emph{finite presentation} if it is the quotient 
of a free algebra on finite generators by an ideal finitely
generated.

{ \bf Free mv-algebras.} 
By standard universal algebra, the free mv-algebra on $n$ generators,
denoted $F_n = F[x_1, \,  
  \ldots \, x_n]$,  is a quotient of the set of  
\emph{terms} in $n$ variables $x_1, \, 
  \ldots \, x_n$ (\cite{COM} section 1.4). 
Two terms $f$, $g$ are considered
equal in $F_n$  if the equation $f = g$ follows from the
axioms. 

\vspace{1ex}

The universal property
characterizing the free algebra means that a term $f \in F_n$ determines by substitution and evaluation a function $A^{n}
\mr{f_A} A$, for any mv-algebra $A$.  A $n$-tuple
$(a_1,\, \ldots \, a_n)$ determines uniquely a morphism $F_n
\mr{\varphi} A$ by defining
  \mbox{$\varphi(f) = f_A(a_1, \ldots \, a_n)$.} This defines an inverse function to the the assignment $\varphi \mapsto
(\varphi(x_1),\, \ldots \, \varphi(x_n))$. There is a bijection
\begin{equation} \label{bijectionell}
\ell:\,[F_n,\, A] \mr{\cong} A^n, \;\;\;\;\; \ell^{-1}(a_1, \ldots \, a_n)(f) = f_A(a_1, \ldots \, a_n).
\end{equation}

\vspace{1ex}

The functions of the form $f_A$
are called \emph{term-functions on $A$}, and the assignment $f
\mapsto f_A$ is a 
morphism $F_n \mr{\delta} A^{A^n}$, where
the exponential notation
stands for the \mbox{mv-algebra} of all functions $A^n \to A$ with the
pointwise structure. The term functions corresponding to the variables
are the projections $\delta{x_i} = {x_i}_A = \pi_i :\,A^n \mr{}
A$. The fact that an 
equation $f 
= g$ holds in a 
mv-algebra $A$ means that the term functions $f_A$ and $g_A$
are equal. 
Clearly, for $A = F_n$, $f = f_{F_n}(x_1,\, x_2,\, \ldots \,
x_n)$, and $F_n$ is isomorphic to the mv-algebra of term-functions on
$F_n$. When the morphism $F_n \mr{\delta} A^{A^n}$ is injective, it establishes
an isomorphism between the free algebra and the algebra of
term-functions on $A$. Thus,\emph{the injectivity of the morphism
$\,\delta\,$ means a
completeness theorem with respect to the algebra $A$.} When it holds it is common usage to say that $A$ \emph{generates} the variety of mv-algebras.

 In particular,
this is the case for $A = [0,\,1]$, and it 
is known as \emph{Chang's completeness
theorem} \mbox{(\cite{COM} 2.5.3).} \emph{We shall abuse notation and write \mbox{$f_{[0,\,1]} = f$.}} The projections and the primitive operations of the mv-algebra
$[0,\,1]$ are continuous (for the usual topology), so all term functions on $[0,\,1]$ are
continuous functions. Thus we have an injective morphism  
\begin{equation}\label{freeisfunctions}
F_n \mr{\delta} Cont([0,\,1]^n,\,
[0,\,1]) \subset [0,\,1]^{[0,\,1]^n}, 
\end{equation}
where continuity is with respect to the usual topology.

Notice that the set $[F_n, \, [0, \, 1]\,] $ is equal by definition to the spectral space  $X_{F_n}$.  Thus, 
$\ell$ establishes a bijection $\ell: X_{F_n} \mr{\cong} [0,\,1]^n$. Under this bijection the evaluations $\widehat{x}_i :X_{F_n} \mr{}  [0,\,1]$ correspond to the projections  $\pi_i: [0,\,1]^n \mr{} [0,\,1]$ ($\pi_i \circ \ell = \widehat{x_i},\; i =
1,\,2,\,\ldots\,n$) so that $\ell$ is continuous for the usual topology in $[0,\,1]$. Since both spaces are compact Hausdorff, we have:
\begin{remark} \label{Xfreeis01}
The bijection $\ell: X_{F_n} \mr{\cong} [0,\,1]^n$ establishes an homeomorphism of topological spaces. Given $f \in F_n$ and a point $p \in  [0,\,1]^n$, $\ell^{-1}(p)(f) = f(p)$.  Under $\ell$  the evaluation function $X_{F_n} \mr{\widehat{f}} [0,\,1]$ corresponds to the term function \mbox{$[0,\,1]^n \mr{f} [0,\,1]$.}
\qed\end{remark}

{\bf Finitely presented mv-algebras.} \label{fpmv}
Since all finitely generated ideals are principal (\ref{mvfacts2} (\ref{aplusbequalab})), a finitely presented mv-algebra is always of the form $R = F_n/(f)$, for some $f \in F_n$.
 A morphism $\varphi$ as in \ref{bijectionell} factors through the quotient if and only if $ f_A(a_1,\,a_2,\, \ldots \, a_n) = 0$. Thus, we have a commutative diagram:
\begin{equation} \label{bijectionellfp}
\xymatrix
              {
                \ell: [F_n,\, A]  \ar[r]^{\cong}  
               & A^n
               \\
                \ell: [R,\, A]  \ar[r]^{\cong} \ar@{}[u] |{\bigcup}
               & Z_{f_A} \ar@{}[u] |{\bigcup}
              }
\end{equation}
where $Z(f_A) \subset A^n$ is the set of zeroes of the term function $f_A$.

\vspace{1ex}

When $A = [0,\,1]$,  we have $Z(f)\subset [0,\,1]^n$
\begin{remark} \label{XfpisZero}
The bijection $\ell$ in \ref{Xfreeis01} restricts to an homeomorphism of topological spaces $\ell: X_{F_n/(f)} \mr{\cong} Z(f)$.
\qed\end{remark}

The injective  arrow $\,\delta\,$ in (\ref{freeisfunctions}) factorizes as follows: 
\begin{equation} \label{fpisfunctions}
\xymatrix
              {
                F_n \ar[r]^(.25){\delta}  \ar[d] 
               & Cont([0,\,1]^n,\, [0,\,1]) \ar[d]
               \\
                R  \ar[r]^(.25){\overline{\delta}} 
               & Cont(Z(f),\, [0,\,1])
              }
\end{equation}
A key non trivial result here is that the arrow $\,\overline{\delta}\,$ is also injective.
\begin{proposition}[Lemma 3.4.8 of \cite{COM}] \label{key}
Given any $f,\,g,\,h \in F_n$:
$$
[g]_{(f)} = [h]_{(f)}  \; \iff  \;  g|_{Z(f)} = h|_{Z(f)},\;\; that  \; is \;\; F_n/(f)\cong F_n|_{Z(f)} .
$$
\end{proposition}
\begin{proof}
The statement is equivalent to  lemma 3.4.8 of \cite{COM} which says:
$$
g \in (f) \;\iff\; Z(f) \subset Z(g).
$$
In fact, it reduces to the lemma when $h = 0$, and it follows from the lemma applied to the function $d(g,\,h)$, where $d$ is the distance operation (\ref{mvfacts1} (\ref{distance}), see also \ref{mvfacts2} (\ref{distance2})).
\end{proof}
\begin{remark} [compare \cite{COM} theorem 3.6.9]\label{fpissemisimple}
In view of remark \ref{XfpisZero} the injectivity of $\overline{\delta}$ amounts to the semisimplicity of finitely presented mv-algebras. Thus finitely presented mv-algebras are strongly semisimple (definition \ref{ssemisimple}).
\end{remark}

We can safely assume the following:
\begin{convention}\label{freeisfunctionsc} $ $
 
 1) We consider the free algebra  \mbox{$F_n = F[x_1, \, x_2,\, 
  \ldots \, x_n]$} to be the algebra of term-functions on $[0,\,1]$. We have then
  \mbox{$F_n \subset Cont([0,\,1]^n,\, [0,\,1])$.}
 
 \vspace{1ex}
  
  2) We consider a finitely presented algebra $R = F_n/(f)$  to be the algebra of term-functions on $[0,\,1]$ restricted to the subset $Z(f) \subset [0,\,1]^n$. We have then
  $R  = F_n|_{Z(f)} \subset Cont(Z(f),\, [0,\,1])$.
\end{convention}

\vspace{1ex}

{\sc The pushout-pullback lemma  for finitely presented algebras.}

Given any finitely presented mv-algebra $R = F_n/(f)$, and $g \in F_n$, we have $R/([g]_f) = (F_n/(f))/([g]_f) = F_n/(f,\,g) = F_n/(f\vee g)$. It follows that the pushout-pullback lemma \ref{pushpull} for finitely presented \mbox{mv-algebras} can be stated as follows:  Given $f_1,\,f_2 \in F_n$, the following diagram is a pullback.
$$
\xymatrix
         {
           F_n/( f_1 \wedge  f_2) \ar[r] \ar[d]
         & F_n/(f_1) \ar[d] 
         \\
           F_n/(f_2) \ar[r]
         & F_n/(f_1 \vee f_2)
         }
$$
By the conventions above and considering that 
$Z(f_1 \wedge  f_2)  = Z(f_1) \cup  Z(f_2)$, and
$Z(f_1 \vee  f_2)  = Z(f_1) \cap  Z(f_2)$, 
 this is then equivalent to the following:
\begin{lemma} \label{pushpullfp}
Given $f_1,\,f_2 \in F_n$, the following diagram of $[0,\,1]$ valued functions is a pullback:
$$
\xymatrix
         {
          F_n|_{Z(f_1) \cup  Z(f_2)} \ar[r] \ar[d]
         & F_n|_{Z(f_1)} \ar[d]
         \\
          F_n|_{Z(f_2)} \ar[r]
         & F_n|_{Z(f_1) \cap  Z(f_2)} 
         }
$$
That is, given $g_1,\,g_2 \in F_n$ such that $g_1|_{Z(f_1) \cap  Z(f_2)} = g_2|_{Z(f_1) \cap  Z(f_2)}$, there exist $g \in F_n$ (necessarily unique upon restriction to $Z(f_1) \cup  Z(f_2)$), such that 
$g|_{Z(f_1)} = g_1|_{Z(f_1)}$ and $g|_{Z(f_2)} = g_2|_{Z(f_2)}$. 
\end{lemma}
\begin{proof}
As in case of hyperarchimedean algebras (lemma \ref{pushpullhyperarchimedean}) we are dealing with algebras of $[0,\,1]$ valued functions. We now use the formula in remark \ref{putaputa}. Set
$$
h = (g_1 \odot  \neg\, nf_1) \vee (g_2 \odot  \neg\, nf_2)
$$
Given $x$,  take $n \geq 1$ such that $nf_1(x) = (n+1)f_1(x)$, $nf_2(x) = (n+1)f_2(x)$, and check (as in \ref{pushpullhyperarchimedean}, \ref{putaputa}) that  for $x \in Z(f_1)$, $h(x) = g_1(x)$, and for $x \in Z(f_2)$, $h(x) = g_2(x)$. The problem now is that we do not have a single $n$ that works for all the $x $ in $Z(f_1) \cup  Z(f_2)$. We eliminate the dependence of the choice of $n$ on $x$ as follows:

\vspace{1ex}

Assume $x \in Z(f_1)$. Then 
$h(x) = g_1(x)  \vee (g_2(x) \odot  \neg\, nf_2(x))$. We shall see there is a $n \geq 1$ such that $g_2(x) \odot  \neg\, nf_2(x) \leq g_1(x)$ for all $x \in Z(f_1)$. Since $g_2(x) \odot  \neg\, nf_2(x) = max \{0,\; g_2(x) - nf_2(x)\}$   (see \ref{mvfacts1} (\ref{[01]})), and $g_1(x) \geq 0$, it follows that it is equivalent to prove
$g_2(x) - nf_2(x) \leq g_1(x)$.

Let $\{T_1,\, \ldots\ T_m\}$ be a set of convex polyhedra as in proposition \ref{partition} below, for \mbox{$f = f_1$,} and  $H = \{f_2,\, g_1,\,g_2\}$. Let $x_{i0}, \, \ldots \, x_{in_i}$ be the vertices of the polyhedron $T_i$. For each $1 \leq i \leq m$ and $0 \leq j \leq n_i$ there is an integer $n_{ij}$ such that $g_2(x_{ij}) - n_{ij}f_2(x_{ij}) \leq g_1(x_{ij})$ for $x_{ij}$. In fact, if $x \in Z(f_2)$, then $g_1(x_{ij}) = g_2(x_{ij})$ and any number $n_{ij}$ will do. If  $f_2(x_{ij}) > 0$, then the inequality will hold if we take a sufficiently large $n_{ij}$. Let $n$ be such that $n_{ij} \leq n$ for all $i, j$. Then $g_2(x_{ij}) - nf_2(x_{ij}) \leq g_1(x_{ij})$ for all $i, j$. Since each $x \in  Z(f_1)$ is a convex combination of the vertices of $T_i$ for some $\, i\,$, and since  $g_1$ and  the function $g_2 - nf_2$ are linear over $T_i$, we get  $g_2(x) - nf_2(x) \leq g_1(x)$. 

For $x \in Z(f_2)$ we do in the same way.
\end{proof}
We state now the key result utilized in the proof of the lemma. This result is  related to the proof of  Lemma 3.4.8 of \cite{COM}. It can be proved as in \cite{COM} proposition 3.3.1, see also \mbox{\cite{NPM} Lemmas 5.2 and 5.3.}
\begin{proposition} \label{partition}
Given any $f \in F_n$ and a finite set $H \subset F_n$, there exists a set of convex polyhedra $\{T_1,\, \ldots\ T_m\}$  whose union coincides with $ Z(f)$,  and such that all the functions $h \in H$ are linear over each $T_i$.
\qed\end{proposition}

\vspace{1ex}

%%%%%%%%%%%%%%%%%%%%%%%%%%%%%%%%%%%%%%%%%%%%%%%%%%%%%%%%%%%%%%%%%%%%%%%
\section{The pushout-pullback lemma} \label{pushpulls}
%%%%%%%%%%%%%%%%%%%%%%%%%%%%%%%%%%%%%%%%%%%%%%%%%%%%%%%%%%%%%%%%%%%%%%%

The general pushout-pullback lemma \ref{pushpull} follows from the
particular case of finitely presented mv-algebras
\ref{pushpullfp}. This is so by categorical nonsense because finite
limits commute with filtered colimits in the category of
mv-algebras. However we find that in this paper a sketch of an
explicit proof is convenient.

Any mv-algebra $B$ is a filtered colimit of finitely presented
mv-algebras. Explicitly, the diagram of all morphisms $R \mr{\alpha}
B$, for all 
finitely presented mv-algebras $R$, is a filtered colimit diagram (with
transition morphisms $(R,\, \alpha) \to (S,\,\beta)$ all $R \mr{\varphi} S$ 
such that $\beta \circ \varphi = \alpha$). Moreover, given \mbox{$a_1,\,a_2 \in B$,} a diagram of shape
 $B/(a_1) \ml{} B \mr{} B/(a_2)$ is in a similar way a filtered
colimit of diagrams  $R/(r_1) \ml{} R \mr{} R/(r_2)$,  $r_1,\,r_2 \in R$, of finitely
presented algebras. It follows that taking the pushouts the resulting squares
also conform a filtered colimit of squares. With this in mind we
proceed to prove the lemma.  
\begin{lemma}[pushout-pullback lemma] \label{pushpull2}
Given any mv-algebra $A$, and two elements $a_1,\,a_2 \in A$, the following pushout diagram is also a pullback diagram (recall that $(a_1,\,a_2) = (a_1 \vee a_2)$). 
$$
\xymatrix
         {
           A/(a_1 \wedge a_2) \ar[r] \ar[d]
         & A/(a_1) \ar[d] 
         \\
           A/(a_2) \ar[r]
         & A/(a_1,\, a_2)
         }
$$
\end{lemma}
\begin{proof}
Let
$b_1, \, b_2 \in A$ be such that they become equal in  $A/ ( a_1
,\, a_2 )$. We have to show there is an element $c \in A$, unique modulo $(a_1 \wedge a_2)$, 
such that $c \mapsto b_1$ in   $A/( a_1 )$, and  
 $c \mapsto b_2$ in  $A/( a_2 )$.

Let $F = F[x_1,\,x_2,\,y_1,\,y_2]$ be the free mv-algebra on four
generators, and consider the morphism $F \mr{} A/(a_1 \wedge a_2)$ determined by the
assignments \mbox{$x_1 \mapsto a_1,\, x_2 \mapsto a_2,\, y_1 \mapsto b_1,\,
y_2 \mapsto b_2$}. This morphism induces the four vertical arrows in
the diagram below. 
$$
\xymatrix@C=4ex@R=2ex
        {
        & F/{( x_1 )} \ar[rr] \ar[dd]
       && F/( x_1,\,x_2 ) \ar[dd] 
       \\ 
          F/( x_1 \wedge x_2 ) \ar[rr]
          \ar[dd]  \ar[ru] 
       && F/( x_2 ) \ar[ru] \ar[dd] 
        & 
       \\ 
        & A/( a_1 ) \ar[rr] 
       && A/ ( a_1 ,\, a_2 ) 
       \\
          A/(a_1 \wedge a_2) \ar[rr] \ar[ru]
       && A/( a_2 ) \ar[ru] 
        & 
        }
$$
The upper square is a pullback, but  the elements $y_1,\,y_2$ do not
become equal in   
$F/{( x_1,\, x_2 )}$. However, they do downstairs in 
$A/ ( a_1 ,\, a_2 )$. Consider the following diagram
(where the outer cube is the cube above):
%\label{colimfiltrante}
$$
\xymatrix@C=4ex@R=2ex
        {
        & F/{( x_1 )} \ar[rr] \ar@{-->}[dd]
       && F/( x_1,\,x_2  ) \ar@{-->}[dd] 
       \\ 
          F/( x_1 \wedge x_2 ) \ar[rr]
          \ar@{-->}[dd]  \ar[ru] %\ar@/_2pc/[dddd]
       && F/( x_2 ) \ar[ru] \ar@{-->}[dd] 
        & 
       \\ 
        & R/( r_1 ) \ar[rr] \ar[dd]|\hole
       && R/( r_1,\, r_2 ) \ar[dd] 
       \\ 
          R/(r_1 \wedge r_2) \ar[rr] \ar[dd]  \ar[ru]
       && R/( r_2 ) \ar[ru] \ar[dd] 
        & 
       \\ 
        & A/( a_1 ) \ar[rr] 
       && A/ ( a_1 ,\, a_2 ) 
       \\
          A/(a_1 \wedge a_2) \ar[rr] \ar[ru]
       && A/( a_2 ) \ar[ru] 
        & 
        }
$$
The square in the bottom is a filtered colimit of the squares (with
vertices finitely presented) in the middle. It follows there is such 
  an square where $y_1,\, y_2$ become already equal in  $R/( r_1,\, r_2
  )$. Since this square is a pullback, there exists an
  element $s \in R$, unique modulo $(r_1 \wedge r_2)$, such that $s \mapsto s_1$ in $R/( r_1
  )$, and $s \mapsto s_2$ in $R/( r_2 )$, where $s_1, \,
  s_2$ are the images by the vertical arrows of $y_1,\, y_2$
  respectively,  $y_1 \mapsto s_1 $, $y_2 \mapsto s_2$. The element $c
  \in A$ image of $s$, $s \mapsto c$, is the required element (recall that the uniqueness is modulo $a_1 \wedge a_2$).   
\end{proof}

%\vspace{1ex}

%%%%%%%%%%%%%%%%%%%%%%%%%%%%%%%%%%%%%%%%%%%%%%%%%%%%%%%%%%%%%%%%%%%%%%%
\section{The compacity lemma.}  \label{posetcoZariski}
%%%%%%%%%%%%%%%%%%%%%%%%%%%%%%%%%%%%%%%%%%%%%%%%%%%%%%%%%%%%%%%%%%%%%%%

The compacity of the prime spectrum space $Z_A$ follows as a consequence of a standard application of Zorn's Lemma. For this we use  a construction of its lattice of open sets which yields a compact locale (see definition \ref{locale}), and then  the space $Z_A$ as the space of points of this locale. This space will be be a compact space provided the locale has enough points (here is when it comes into play Zorn's Lemma).

\vspace{1ex}

\begin{center}{\sc Construction of the prime spectrum of
  mv-algebras via sheaf theory of  posets}\end{center}
  
  \vspace{1ex}

Any mv-algebra $A$ is a distributive lattice. On the other hand, the
principal ideals under subset inclusion also form a distributive
lattice which is a quotient lattice of $A$.  In fact, it is the quotient
lattice determined by the following equivalence relation:
\begin{definition} \label{sim}
Given any mv-algebra $A$ and two elements $a,\, b \in A$, 
$$
a \sim b \;\iff\; (a) = (b) \;\iff\; \exists\,n \;\;|\;\; a \leq nb
\;\;and\;\; b \leq na
$$
\end{definition}  
\begin{proposition} \label{supequaloplus}
The relation defined in \ref{sim} is a lattice congruence:
$$
(a_1) = (a_2),\;(b_1) = (b_2) \;\;\rimply \;\; (a_1 \wedge b_1) =
(a_2 \wedge b_2),\;\;  (a_1 \vee b_1) = (a_2 \vee b_2).  
$$
\qed\end{proposition}
This proposition immediately follows from the following:
\begin{proposition}
Given any mv-algebra $A$, and two elements $x,\,y \in A$:
$$
(x \wedge y) = (x) \cap (y), \;\;\;\;\; (x \vee y) = (x,\,y) = (x
\oplus y).
$$
\end{proposition}
\begin{proof}
The first equation follows from \ref{mvfacts3} (\ref{p2}), and the
second since \mbox{$x \wedge y = 0 \;\iff \; x = 0,\; y = 0 \;\iff \;
  x \oplus y = 0$.} 
\end{proof}

\vspace{1ex}

Given any distributive lattice, the opposite order also defines a
distributive lattice. We shall denote by $A^{op}$ the
distributive lattice determined by the opposite order in a mv-algebra
$A$. We consider the opposite lattice of the lattice of principal
ideals defined above:
\begin{definition}
Given any mv-algebra $A$, we denote by $V_A$ the quotient of the
lattice $A^{op}$ by the congruence defined in \ref{sim}, $A^{op} \mr{}
V_A$. 

The quotient map is denoted by an over-lining, $a \mapsto
\overline{a}$. In this way, 
$$a \leq b \;\rimply \overline{b} \leq
\overline{a},\;\;\;\overline{a \wedge b} = \overline{a} \vee
\overline{b},\;\;\; and \;\;\;\overline{a \vee b} = \overline{a} \wedge
\overline{b}.
$$
\end{definition}

We refer to the definition of point of a inf-lattice (\ref{latticepoints}).
\begin{proposition} \label{pointsofVA}
A point of the lattice $A^{op}$,  $A^{op} \mr{p} 2$ satisfies the
equation 
$$
p(a \oplus b) = p(a) \wedge p(b)
$$ 
if and only if it
factorizes by the quotient as indicated in the diagram:
$$
\xymatrix
         {
           A^{op} \ar[r] \ar[rd] 
         & V_A \ar@{-->}[d] 
         \\
         & 2
         }
$$
\end{proposition}
\begin{proof}
 Assume that $p$ factorizes, then (abusing notation) we can write $p =
 \overline{p}$. From proposition \ref{supequaloplus} it follows
 $\overline{a\oplus b} = \overline{a \vee b} = \overline{a} \wedge
 \overline{b}$. Then, \mbox{$p(a 
 \oplus b) = p(\overline{a \oplus b}) = p(\overline{a}) \wedge
 p(\overline{b}) = p(a) \wedge p(b)$.}
Assume now the equation. It follows that $p(x) = p(nx)$. We have to
see that if $a \sim b$, then $p(a) = p(b)$. But under the hypothesis,
$b \leq na$ and $a \leq mb$. Thus $p(a) = p(na) \leq p(b)$, and   $p(b)
= p(mb) \leq p(a)$.    
\end{proof}

We consider the lattices $A^{op}$ and $V_A$ furnished with the
Grothendiek topology $\jmath_f$ of finite suprema (see \ref{gtopology},
\ref{lattice}). Then, proposition \ref{pointsofVA} proves the following.
\begin{proposition} \label{pointsofVA2}
The points of the site $(V_A,\, \jmath_f)$ are exactly the
prime ideals of the mv-algebra $A$. 
\qed\end{proposition}
From this and theorem \ref{isopoints} (\ref{localiso}) it follows
\begin{proposition}
The topological space $Z_A$ (the base of the mv-space
$\ff{Spec}_A$, see \ref{Spec_A}) is the same that the space
$P_\jmath(V_A)$ of points  of the site $(V_A,\, \jmath _f)$.
\qed\end{proposition}
Then, theorem \ref{spectralspace} proves, in particular, the compacity
lemma \ref{compacity}.
\begin{corollary} \label{Z_Acompact}
Given any mv-algebra $A$, the spectral space $Z_A$ is sober, compact,
and has a base of compact open sets.  
\qed\end{corollary}
We  finish this section characterizing the open sets of this space. We
know that they correspond to the elements of the locale $I_\jmath(V_A)$, that
is, the lattice ideals of $V_A$ (see \ref{lattice}). These in turn are
certain lattice ideals of $A^{op}$, that is, lattice filters of
$A$. Considering that given any $x \in 
 A$, $(x) = (nx)$, 
we let the reader prove the following
characterization.
\begin{proposition}
The open sets of the spectral space $Z_A$ correspond to the lattice
filters $U \subset A$ which have the descent property \mbox{``$\,na \in U
\;\rimply \; a \in U$''.} The corresponding open set is given by $U =
\{P\;|\; \exists\, a \in U,\; a \in P\}$.
\end{proposition}   

%%%%%%%%%%%%%%%%%%%%%%%%%%%%%%%%%%%%%%%%%%%%%%%%%%%%%%%%%%%%%%%%%%%%%%
\section{The theorem of McNaughton} \label{McN}
 %%%%%%%%%%%%%%%%%%%%%%%%%%%%%%%%%%%%%%%%%%%%%%%%%%%%%%%%%%%%%%%%%%%%%%

Recall that we consider the free mv-algebra $F_n = F[x_1, \, x_2,\, 
  \ldots \, x_n]$ to be the algebra of term functions $F_n \subset
[0,\,1]^{[0,\,1]^n}$ (\ref{freeisfunctionsc}). Recall the following:
\begin{proposition}[\cite{COM} 3.1.9] \label{sharp}
A linear function with integer coefficients \mbox{$h = s_0 + s_ 1x_1 +
\ldots + s_nx_n$, $s_i \in  \bb{Z}$,}  determines a term function, denoted $h^\sharp$, by means of the definition    
$h^\sharp = (h \vee 0) \wedge 1$.
\qed\end{proposition}

Any convex polyhedron  $P \subset [0, \, 1]^n$ is the intersection of $[0, \, 1]^n$ and a finite set  of closed half spaces defined by linear functions. It follows:
\begin{remark} \label{zeroset}
A convex polyhedron $P \subset [0, \, 1]^n$  with rational vertices is the Zero set of a term function $f = h_1^\sharp \,\vee \ldots \, \vee \, h_k^\sharp$, 
$P = Z(f)$, where $h_1, \, \ldots \, h_k$ are linear functions with with integer coefficients.
\qed\end{remark}
As a particular case of proposition \ref{partition} we have
\begin{proposition}
Given any $f \in F_n$, there exists a set of convex polyhedra $\{T_1,\, \ldots\ T_m\}$  whose union coincides with $ Z(f)$.
\end{proposition}

Remark that given a finite set of term functions $f_1,\, \ldots \, f_m \, \in \, F_n$, the open sets $W_{f_i} = Z(\widehat{f}_i) \subset Z_{F_n}$ cover the prime spectral space $Z_{F_n}$ exactly when $f_1 \wedge \ldots \wedge f_m = 0$, which in turn is equivalent to the fact that the zero sets of the term functions $Z(f_i) \subset [0, \, 1]^n$ cover the cube $[0, \, 1]^n$.

By refining the covers if necessary, we have
\begin{remark} \label{open=poly}
The (finite) open covers of the prime spectral space $Z_{F_n}$ correspond to the (finite) covers of the cube $[0, \, 1]^n$ by convex polyhedra with rational vertices. 
\end{remark}

\vspace{1ex}

Recall that a McNaughton function with constituents in a finite set $\{h_1, \, \ldots \, h_k\}$ of linear functions with integer coefficients is a continuous function $[0,\,1]^n \mr{\tau} [0,\,1]$ such that for each point $x \in [0,\,1]^n$, $\tau(x) = h_i(x)$ for some $i \in \{1,\ldots k\}$ (see \cite{COM} definition 3.1.6). The following is not difficult to prove (compare with proposition
\ref{partition}).
\begin{proposition}[\cite{COM} 3.3.1]\label{partition0}
Given a McNaughton function $\tau$ with linear constituents 
$h_1, \, \ldots \, h_k$, there exists a set  of convex polyhedra  with rational vertices $\{T_1,\, \ldots\ T_m\}$ whose union coincides with $[0, \,1]^n$, and such that for each $T_i$ there is a $h_j$ with $(\tau = h_j)|_{T_i}$.
\end{proposition}
In conclusion  we have:
\begin{proposition} \label{partition2}
Any McNaughton function $\tau$ is determined by a cover of $[0, \, 1]^n$ by convex polyhedra $T_i = Z(f_i)$, $f_1,\, \ldots \, f_m \, \in \, F_n$, $f_1 \wedge \ldots \wedge f_m = 0$, and a compatible family  $g_1,\, \ldots \, g_m \in F_n$ of term functions, $(g_i = g_j)|_{T_i \cap T_j}$. Vice-versa, any such data determines a McNaughton function by setting  $(\tau = g_i)|_{T_i}$ (the second assertion is justified by \ref{freeisMc} below).
\end{proposition}

Clearly the projections $x_i$ and the constant function $0$ are
McNaughton functions, and directly from the definition it can be easily seen that  \mbox{McNaughton}
functions form a \mbox{mv-subalgebra} of $[0,\,1]^n$. It follows:
\begin{proposition}[\cite{COM} 3.1.8] \label{freeisMc}
Term functions are McNaughton functions. That is, $F_n$ is a subalgebra of $M_n$, $F_n \subset M_n$ (where $M_n$ denotes the \mbox{mv-algebra} of McNaughton functions)
\qed\end{proposition}
\emph{McNaughton's Theorem establishes the converse result}, that is,
that every McNaughton function is a term function, $F_n \supset M_n$.

\vspace{1ex}
 
It is convenient now to write explicitly the definition of global sections of the prime spectrum (see \ref{sections}). 
\begin{fact} \label{sectionsoffree}
A global section of the prime spectrum $\sigma \in \Gamma(Z_{F_n},\, E_{F_n})$ is determined by a cover of $Z_{F_n}$ by open sets $W_{f_i} = Z(\widehat{f}_i)$, $f_1,\, \ldots \, f_m \, \in \, F_n$, $f_1 \wedge \ldots \wedge f_m = 0$, and a compatible family  $\,g_1,\, \ldots \, g_m \in F_n$ of term functions, $(\widehat{g}_i = \widehat{g}_j)|_{W_{f_i} \cap W_{f_j}}$. In such case, $(\sigma = \widehat{g}_i)|_{W_{f_i}}$.
\qed\end{fact}

We see  that a datum $((f_1,\, \ldots \, f_m),\,  (g_1,\, \ldots \, g_m))$, $f_1 \wedge \ldots \wedge f_m = 0$,  \mbox{determines} either a McNaughton function or a global section, depending on whether \mbox{$(g_i = g_j)|_{Z(f_i) \cap Z(f_j)}$} or 
$(\widehat{g}_i = \widehat{g}_j)|_{Z(\widehat{f}_i ) \cap Z(\widehat{f}_j )}$.
But  \cite{COM} Lemma 3.4.8 says precisely that this two conditions are equivalent (see \ref{compatible} and \ref{key}). This shows on the spot the equality (recall $ (Z_{F_n}, E_{F_n}) = \ff{Spec}_{F_n}$)
$$M_n = \Gamma(\ff{Spec}_{F_n})$$

We establish now a precise statement of this fact.  Consider the composite morphism (see \ref{kappalambda} and \ref{Xfreeis01})
$$
\Gamma(Z_{F_n}, E_{F_n}) \mr{\kappa^*\lambda_*}  Cont(X_{F_n}, \, [0,\,1]) \mr{(\ell^{-1})^*}  Cont([0,\,1]^n,\, [0,\,1])
$$
This morphism sends a global section $\sigma$ into the function $\tau$ defined by 
$$
\tau(p) = \lambda\sigma\kappa (\ell^{-1} (p)), \;\;\; p \in [0,\,1]^n
$$  
In  fact: $(\ell^{-1})^* \kappa^*\lambda_*(\sigma)(p) = \kappa^*\lambda_* (\sigma)(\ell^{-1}(p)) = \lambda\sigma\kappa (\ell^{-1} (p))$.

\vspace{1ex}

Given $g \in F_n$, we have 
$\lambda\,\widehat{g}\,\kappa (\ell^{-1} (p)) = 
\ell^{-1}(p)(g) = g(p)
$.
Thus it sends $\,\widehat{g}\,$ into $g$. It follows it sends the datum associated to a global section into the same datum, but now defining a McNaughton function. This shows it is surjective into the mv-algebra $M_n$. From proposition \ref{strongsemi2} and remark \ref{fpissemisimple} it follows that it is also injective. This completes the proof of the following: 
\begin{theorem} \label{mc=gamma}
The morphism $\Gamma(Z_{F_n}, E_{F_n}) \mr{} Cont([0,\,1]^n,\, [0,\,1])$ which sends a global section $\sigma$ into the function $\tau(p) = \lambda\sigma\kappa (\ell^{-1} (p))$, establishes an isomorphism $\ff{Spec}_{F_n} \mr{\cong} M_n$ between the mv-algebra of global sections of the prime spectrum of the free algebra and the mv-algebra of McNaughton functions. This isomorphism sends a global section of the form $\widehat{g}$ into the term function $g$.
\qed\end{theorem}
If the reader wants to identify $[0,\,1]^n$ to a subset of $Z_{F_n}$, and the fibers of $E_{F_n}$ over a maximal ideal with the interval $[0,\,1]$, then he can think this isomorphism as the restriction of $\sigma$ to $[0,\,1]^n$. In this way, a McNaughton function has a unique extension into the whole prime spectrum $Z_{F_n}$.

The global sections of the maximal spectrum $M_{F_n}$ (also identified
with the cube $[0,\,1]^n$ but this time furnished with the coZariski
topology, which has as a base of open sets the Zero sets of the term functions), like McNaughton functions, are given by linear functions on convex polyhedra (this time a possible infinite family), but unlike McNaughton functions, they are only continuous for the (much) finer coZariski topology. They do not extend to the whole prime spectrum. For example, the function which is equal to $1$ in $\{0\} = [0,\, 0]$, and constantly zero on each interval $[1/(n+1), \, 1/n]$, $n \in \bb{N}$,  is such a function.

\vspace{1ex}

From theorem \ref{mc=gamma} it clearly follows the following corollary
\begin{theorem}
The theorem of McNaughton is equivalent to the representation theorem \ref{main} for free mv-algebras.
\qed\end{theorem}
This shows that the representation theorem can be viewed as a vast
generalization of the McNaughton theorem, from free mv-algebras to
arbritrary mv-algebras. In particular, the representation theorem
furnish a proof of the theorem of McNaughton.

%
%be a set of compact convex $n$-dimensional simplexes, any two $T_i$ and 
%$T_j$ being either disjoint or intersecting in a comun face,
%
%%%%%%%%%%%%%%%%%%%%%%%%%%%%%%%%%%%%%%%%%%%%%%%%%%%%%%%%%%%%%%%%%%%%%%%
\section{Appendix. Sheaf theory of posets} \label{appendix}
%%%%%%%%%%%%%%%%%%%%%%%%%%%%%%%%%%%%%%%%%%%%%%%%%%%%%%%%%%%%%%%%%%%%%%%

In this appendix we fix notation, terminology, and prove some general
results. We need this results 
as they are stated here. They are
known and can be found (in one form 
or another) in the literature. However, the development here is original
in print as far as we know. We did not find any treatment which arrive fast to the results, and that in the statements and the proofs is  
faithful to the idea we want to stress here, namely, that of a \emph{sheaf
theory of posets}. With posets in the role of categories, the two
elements poset in the role of the category of sets,  
inf-lattices in the role of categories with finite limits, and locales
in the role of Grothendieck topoi.

\vspace{1ex} 

We shall consider a partial order to be a reflexive and transitive
relation, not necessarily antisymmetric. A set furnished with such a
relation will be called a \emph{poset}. This is equivalent to a
category taking its homsets in the poset $2 = \{0,\, 1\} =
\{\emptyset,\, \{*\}\}$. As usual, $x \leq y \;\iff\; hom(x,\, y)
\not= \emptyset$. Under this equivalence a functor is the same thing
that an order preserving function.  Given any poset $H$, under the
usual bijection between 
subsets and characteristic functions, functors $H
\mr{p} 2$ correspond with poset-filters \mbox{$P \subset H$.} Functors
$H^{op} \mr{u} 2$ are called \emph{presheaves}, and correspond with
poset-ideals $U \subset H$.  The set $I(H)$ of all ideals, ordered by
  inclusion, \mbox{$I(H) = 2^{H^{op}}$,} is a
  \emph{locale} (see \ref{locale} below). The locale structure is
  given by the \emph{union} and \emph{intersection} of subsets. There is a 
\emph{Yoneda} functor $H \mr{h} I(H)$, sending an element $a \in
H$ to the principal ideal $(a]$. This functor is \emph{full}, meaning
  that for any $x, \, y \in H$, we have $x
  \leq y \;\iff\; h(x) \leq h(y)$.

\vspace{1ex}
 
Following Joyal-Tierney \cite{JT}, we  think of locales as dual
objects for generalized (may be pointless) topological spaces, the
local being its 
lattice of open sets. In the same vein, we think of
\emph{inf-lattices} as open basis for locales. 

\pagebreak

Recall:
\begin{definition} \label{locale} A \emph{locale} is a complete lattice in which
  finite infima 
  distribute over arbitrary suprema. A morphism of locales $L
  \mr{f^*} R$ is defined as a function $f^*$ preserving finite infima
  and arbitrary suprema (we put an upper star to indicate that these
  arrows are to be considered as inverse images of a map between the
  formal duals $\overline{R} \mr{f} \overline{L}$. The formal dual of
  a local is called an \emph{space} in \cite{JT}, we shall call it a
  \emph{localic space}.
\end{definition}
\begin{remark} \label{opensubspace}
Given a local $L$, each element $u \in L$ determines a local $L_u =
\{x \,|\, x \leq u\}$. Notice that the inclusion $L_u \subset L$ is
not a morphism of locales since does not preserve $1$. There is a
quotient morphism of locales $L \to L_u$ given by $x \mapsto x \wedge
u$. This determines the open subspace $\overline{L}_u \hookrightarrow
\overline{L}$. 
\end{remark}

A topological space $X$ is considered
  to be a localic space by means 
  of its lattice  $\cc{O}(X)$ of open sets. A continuous function $Y
  \mr{f} X$ determines a morphism of locales
  $\cc{O}(X) \mr{f^*} \cc{O}(Y)$ given by the inverse image.

The poset $2 = \{0,\,1\} =
  \cc{O}(\{*\})$ is the 
  \emph{singleton} or \emph{terminal} localic space. Given any locale
  $L$, there exists a
  unique locale morphism  $\{0,\,1\} \to L$, $1 = \overline{2}$.

\begin{definition} A \emph{point} $p$ of a locale $L$ is a
  locale morphism $L \mr{p^*} 2$ 
  \mbox{(that is, $1 \mr{p} \overline{L}$).} 
\end{definition}
\begin{proposition} \label{spaceofpoints}
Given any local $L$, the set of
  points $P_L$ has a canonical topology whose open sets are the
  subsets $W_u \subset P_L$, $W_u = \{p\,|\, p^*{u} = 1\}$, for
  $u \in L$. 
\qed\end{proposition}
There is a surjective morphism of locales $L \mr{\rho}
  \cc{O}(P_L)$. 
\begin{definition} \label{epoints} We say that a 
  local $L$ has \emph{sufficiently many (or enough)} points when $\,\rho\,$ 
  is injective, $u \neq v \;\rimply\; W_u \neq W_v$.  That is:
$$
u \neq v \;\rimply\; \exists \, p \,
    |\, p^*u = 1,\; p^*v = 0. 
$$
In this case, the localic space $\overline{L}$ is topological, $L
\mr{\cong} \cc{O}(P_L)$.
\end{definition}
The topological space $P_L$ is a
  \emph{sober} space (that 
  is, every nonempty irreducible closed subset has a unique generic
  point, \cite{J} chapter II). The category of sober topological spaces is dual to the
  category of locales with enough points. 

\vspace{1ex}
      
A inf-lattice is a poset with finite
infima (in particular, the 
empty infimun or top element 1). A morphism of inf-lattices is a
infima  preserving (thus also order preserving) function.

\begin{definition} \label{latticepoints}
A \emph{point} of a inf-lattice $V$ is a morphism $V \mr{p} 2$. A
\emph{presheaf} is a order reversing  map $V^{op} \mr{u}
2$. Points correspond to inf-lattice filters, and presheaves to poset
(not necessarily inf-lattice) ideals \mbox{(see \ref{2points} and
\ref{2sheaves} below).} 
\end{definition} 
\begin{definition} \label{gtopology}
A \emph{grothendieck topology} $\jmath$ in a inf-lattice $H$ consists
of, for every 
$a \in H$, a set $\jmath(a)$ of families $a_i \leq a$, called \emph{covers},
subject to the following axioms:

\pagebreak 

$ \hspace{5ex} i) \;  x \cong a \; \in \jmath(a). $

 \vspace{1ex}
                          
$ \hspace{5ex} ii) \; a_{i,\,j} \leq a_i \, \in \jmath(a_i),\;  a_i \leq a \, \in
\jmath(a)  \;\;\rimply \;\; a_{i,\,j} \leq a \, \in \jmath(a). $

\vspace{1ex}
   
$ \hspace{5ex} iii)\;  a_i \leq a \, \in \jmath(a),\; b \in H  \;\;
\rimply \;\; a_i 
\wedge b \,\in  \jmath(a \wedge b).$ 

\vspace{1ex}

The topology is said to be \emph{subcanonical} if the covers are
suprema, that is: 

\vspace{1ex}

$\hspace{5ex} iv) \, for\; every \;\; a_i \leq a \in \jmath(a), \; a =
\bigvee_i \, a_i$.

\vspace{1ex}

An inf-lattice furnished with a Grothendieck topology is called a \emph{site}.
\end{definition}
We shall often say \emph{topology} instead of  Grothendieck
topology. There is a minimal or trivial topology whose covers are the
isomorphisms. In general, 
it is possible for some elements $a$ to be covered by the empty
family, we denote this by $\emptyset \in
\jmath(a)$. In practice this is the case when the
inf-lattice has a bottom 
element (the empty supremun) $0 \in H$. Then, usually it is required    
$\emptyset \in \jmath(0)$. 
\begin{definition} \label{2points} Let $(H,\, \jmath)$ be a site. A
  \emph{point} is a inf-preserving functor $H \mr{p} 
  2$ such that it sends covers into epimorphic families. Writing \mbox{$P =
  \{a\,|\, p(a) = 1\}$,} points
  correspond with $\jmath$-prime inf-lattice filters.  These are   
  subsets \mbox{$P \subset H$} such that:
\begin{eqnarray*}
i)    &  1 \in P,                                     \\
ii)   &  a \geq b \in P & \rimply  \; a \in P,        \\
iii)  & a,\, b \in P    & \rimply  \; a \wedge b \in P,    \\
iv.a)  & a_i \leq a \in \jmath(a) \; and \; a \in P & \rimply  \;
\exists \, i \; | \; a_i \in P.                        \\
iv.b)  &  \emptyset \in \jmath(a) & \rimply \; a \notin P.
\end{eqnarray*}
\end{definition}
Remark that such a filter $P$ may not be \emph{proper}. However, in
practice, $\jmath$-prime filters are proper. 
Usually there is a bottom element $0 \in H$, and
$\emptyset \in \jmath(0)$, so that $0 \notin P$ for any $\jmath$-prime
filter $P$.
\begin{definition} \label{2sheaves} Let $(H,\, \jmath)$ be a site. A
  \emph{sheaf} is a presheaf $H^{op} \mr{u} 2$ such that:

\vspace{1ex}

$ \emph{Sheaf axiom:} \; a_i \leq a \, \in \jmath(a) \; and \;\; \forall \,i
\; u(a_i) = 1  
\;\;\; \rimply  \;\;\;u(a) = 1.$

\vspace{1ex}
  
Writing $U = \{a \,|\, u(a) = 1\}$, sheaves correspond with
$\jmath$-ideals. These are poset  (not
inf-lattice) ideals satisfying the sheaf axiom. That is, subsets
\mbox{$U \subset H$} such that:  

\vspace{1ex}
                            
$ i)\,  a \leq b \in U \;\; \rimply \;  a \in U.$

\vspace{1ex}
  
$ ii.a)\,  a_i \leq a \, \in \jmath(a) \; and \;\; \forall \,i \; a_i \in U 
\;\;\; \rimply  \;\;\;a\in U.$

\vspace{1ex}
 
$ii.b)\, \;\emptyset \,\in \jmath(a) \;\;\; \rimply  \;\;\; a \in U$.
\vspace{1ex}
\end{definition}
Usually there is a bottom element $0 \in H$, and
$\emptyset \in \jmath(0)$, so that $0 \in U$ for any $\jmath$-ideal
$U$.

\vspace{1ex}

The reader may examine now the example \ref{lattice} below.

\vspace{1ex}

We consider next the poset  $I_\jmath(H)$ of all $\jmath$-ideals, ordered
by inclusion, $I_\jmath(H)\subset I(H)$. We collect in the following
theorem the basic properties of this 
construction. We leave the reader to check the details of their proof. It is a
instructive but straightforward task.
\begin{theorem} \label{ideallocalefacts} Given any site $(H,\, \jmath)$:

\begin{enumerate}
\item \label{asheaf}
Given any subset $S \subset H$, the set 
$$
\hspace{6ex} \# S = \{x \,|\, \exists \; a_i \leq a \, \in
\jmath(a),\; a_i \in S \; \forall \,i, \;\; x \leq a \}
$$
is a $\jmath$-ideal, said to be \emph{generated} by $S$. 

\vspace{1ex}

\item \label{idealslocale}
The poset $I_\jmath(H)$ is a locale, $I_\jmath(H)\subset I(H)$. The generated
  $\jmath$-ideal determines a morphism of locales $I(H) \mr{\#}
  I_\jmath(H)$, such that for $S \in I(H), \; U \in I_\jmath(H)$: 
$$
\hspace{8ex} \#S \leq U \;\iff\; S \leq U\;\;\;\;(\# \;is\; left\;
adjoint\; to\; 
the\; inclusion). 
$$ 
The locale structure of $I_\jmath(H)$  
  is given by the following: 
  $$U \wedge V = U \cap V, \;\;\;\; \bigvee_i \, U_i \, =
  \#\bigcup_i \, U_i$$
\item \label{topbot}
The bottom element is the $\jmath$-ideal $\#\emptyset = \{a \,|\, \emptyset \in
\jmath(a)\}$, and the top element is the
whole set $H$. Given a $\jmath$-ideal $U$, $U = H \;\iff 1 \in U$.

\vspace{1ex}

\item \label{epsilon}

The composite $H \mr{h} I(H) \mr{\#} I_\jmath(H)$ determines a
\mbox{inf-lattice} morphism $\varepsilon = \# h$ sending covers into
suprema. 
Given $a \in H$, 
$$\hspace{4ex} \varepsilon(a) = \#(a] = \{x
\;|\; \exists \; 
b_i \leq b \, \in 
\jmath(b),\; b_i \leq a\; \forall \,i, \; x \leq b \}.
$$

\vspace{1ex}

\item \label{density}
Given any $U \in I_\jmath(H)$, $U = \bigvee_{a \in U} \,
\varepsilon(a)$. Thus the elements of the form $\varepsilon(a)$, $a \in H$,
are a \emph{base} of the locale $I_\jmath(H)$ (notice that
\mbox{$\varepsilon(a \wedge b) = \varepsilon(a) \wedge \varepsilon(b)$).}

\vspace{1ex} 

\item
The topology is subcanonical if and only if the segment  $(a]$ is
    already a $\jmath$-ideal. That is, $\varepsilon(a) = (a]$. This is the
      case if and only if $\varepsilon$ is full, that is, for any $x,
      \, y \in H$, $x 
  \leq y \;\iff\; \varepsilon(x) \leq \varepsilon(y)$.
\end{enumerate}
\end{theorem}
\begin{proof}
Clearly $I(H) = 2^{H^{op}}$ is a locale (in fact, it has the pointwise
structure determined by the locale $2$). Next, check that the
generated sheaf $\#$ preserves finite infima and that it is left
adjoint to the inclusion. From this it easily follows that
$I_\jmath(H)$ is a locale. The rest is straightforward.
\end{proof}
%The following theorem and remark are esily proved. the reader should
%check first the validity of the remark, and use the remark to check
%the theorem.
\begin{theorem} \label{isopoints}
 Given any site $(H,\, \jmath)$:

\begin{enumerate}
\item \label{universal}
Given any \mbox{locale} $L$ and a inf-lattice-morphism $H \mr {f} L$ sending
covers into suprema, there exists a unique morphism of locales
$I_\jmath(H) \mr{f^*} L$ such that $f^* \varepsilon = f$. 

$f^*$ is
determined by the formula $f^*(U) = \bigvee_{a \in U}\,
f(a)$. 

Furthermore, for any two $f,\;g$, $\;f \leq g \;\iff\; f^* \leq g^*$.

\vspace{1ex}

\item \label{localiso}

In particular composition with $\varepsilon$
      establishes a  bijection

\noindent $P_{I_\jmath(H)} \mr{\cong} P_\jmath(H) =\;$ \{$\,P \subset H\,|\, P
  \;is\;a\; \jmath$-prime inf-lattice filter\}. 

\vspace{1ex}

\noindent The topology of
  $P_{I_\jmath(H)}$  induces a
  topology in the set $P_\jmath(H)$. A base for this topology is given by the
  sets \mbox{$W_a = W_{\varepsilon(a)} = \{P \,|\, a \in P\}$}, for
  each $a
  \in H$ (see proposition \ref{spaceofpoints} and theorem
  \ref{ideallocalefacts} (\ref{density})).  
\end{enumerate}
\end{theorem}
\begin{proof}
We check that $f^*$ preserves finite infima. Any order
preserving map satisfies $f^*(U \cap V) \leq f^*(U) \wedge f^*(V)$. For
the other direction we do as follows:

$f^*(U) \wedge f^*(V) =
\bigvee_{a \in U}\, 
f(a)\; \wedge\;  \bigvee_{a \in V}\,f(a) \;=\;  \bigvee_{a \in U, \, b
  \in V}\,f(a) \wedge f(b) \;=\;\bigvee_{a \in U, \, b
  \in V}\,f(a \wedge b) \;\leq \; f^*(U \cap V) $. The proof of the
other statements in the theorem is straightforward and it is left to the reader.
\end{proof}  
\begin{remark} \label{fU=1}
Let $H \mr{p} 2$ be a point with corresponding $\jmath$-prime
inf-lattice filter $P \subset H$, and $U \subset H$ be any
$\jmath$-ideal. Then (see \ref{2points}, \ref{2sheaves}):

\vspace{1ex}

$
  \hspace{18ex} p^*(U) = 1 \;\;\; \iff \;\;\;U \cap P \neq \emptyset.
$ 
\qed\end{remark}
\begin{proposition} \label{compactlocale}
Given any site $(H,\, \jmath)$ whose covers are \emph{finite
    families}, for any $a \in H$, $\varepsilon(a) \in I_\jmath(H)$ is
  \emph{compact}. That is:
$$
\hspace{6ex} \varepsilon(a) \leq \bigvee_i \, U_i \; \rimply \; \exists \;
  i_1,\, i_2,\; \ldots \; i_n \;\;  | \;\; \varepsilon(a) \leq U_{i_1}
  \vee U_{i_2} \vee \; \ldots \; \vee U_{i_n}.
$$   
Thus, the locales $I_\jmath(H)_{\varepsilon(a)}$ are compact. In
particular (for a = 1), $I_\jmath(H)$ is a compact locale with a base
of compact elements (see \ref{ideallocalefacts} (\ref{density})).
\end{proposition}
\begin{proof}
The following
chain of equivalences (justified by \ref{ideallocalefacts} (\ref{asheaf}),
(\ref{idealslocale}) and  (\ref{epsilon})) proves
the proposition: 
\begin{center}
\begin{picture}(0, 144) 
\put (-40, 132)  {$\varepsilon(a) \leq \bigvee_i \, U_i$} 
\put (-150, 124){\line(1, 0){305}}
\put (-145, 110) {$\#(a] \leq \bigvee_i \, U_i \;\;\equiv\;\;(a] \leq
  \bigvee_i \, U_i \;\;\equiv\;\; a \in \bigvee_i \, U_i \;\;\equiv\;\; a \in
  \# \bigcup_i \, U_i$}
%\put (-160, 56){$a$} 
\put (-150, 102){\line(1, 0){305}}
\put (-138, 88)  {$\exists \, a_{i_1},\;\;\ldots\;\; a_{i_n} \leq
  b \;\in\; \jmath(b),\;\; a_{i_1},\;\;\ldots\;\; a_{i_n} \in
  \bigcup_i \, U_i \,,\;\; a \leq b$}
\put (-160, 80){\line(1, 0){325}}
\put (-155, 66)  {$\exists \, a_{i_1},\;\;\ldots\;\; a_{i_n} \leq
  b \;\in\; \jmath(b),\;\; a_{i_1},\;\;\ldots\;\; a_{i_n} \in
  U_{i_1} \cup \;\ldots\; \cup U_{i_n} \,,\;\; a \leq b$}
\put (-160, 58){\line(1, 0){325}}
\put (-115, 44)  {$a \in  \#(U_{i_1} \cup \;\ldots\; \cup
  U_{i_n})\;\;\equiv\;\; a \in  U_{i_1} \vee \;\ldots\; \vee
  U_{i_n}$}
\put (-120, 36){\line(1, 0){240}}
\put (-118, 22){$(a] \leq  U_{i_1} \vee \;\ldots\; \vee
  U_{i_n} \;\;\equiv\;\; \#(a] \leq  U_{i_1} \vee \;\ldots\; \vee
  U_{i_n}$}
\put (-120, 14){\line(1, 0){240}}
\put (-50, 0){$\varepsilon(a) \leq U_{i_1} \vee \;\ldots\; \vee
  U_{i_n}$}
\end{picture}
\end{center}
\end{proof}
Given a compact locale $L$, the formal dual is (by definition) a
compact localic space. The topological space $P_L$ (see
\ref{spaceofpoints}) need not be a compact topological space unless
$L$ has 
enough points, in which case $L \cong \cc{O}(P_L)$. As usual, the
existence of ``points'' depends of the axiom of choice. The following
theorem follows by an standard choice argument:
\begin{theorem}
Given any site $(H,\, \jmath)$ whose covers are \emph{finite
    families}, the locale $I_\jmath(H)$ has enough points. 
By remark \ref{fU=1}, this amounts to the following statement:

\vspace{1ex}

Given any two j-ideals
   $U, \; V \subset H$, if $U \neq V$, then there exist a
  $j$-prime inf-lattice filter $P \subset H$ such that $U \cap P
  = \emptyset$, and $V \cap P \neq \emptyset$.  
\end{theorem}
\begin{proof}
Take an element $a \in H$ such that  $a \in V$, $a \notin U$. Consider
the set $\cc{F}$ of inf-lattice 
filters $F \subset H$, $\cc{F} = \{F \,|\, a \in F,\; U \cap F = \emptyset
\}$. Clearly, if $F \in \cc{F}$,  $U \cap F
  = \emptyset$, and $V \cap F \neq \emptyset$. We shall see that there
  is a $\jmath$-prime filter in $\cc{F}$. 

\vspace{1ex}

The inf-lattice filter
$[a) \subset H$, $[a) = \{x\,|\, a \leq x\}$ is in $\cc{F}$, so
    $\cc{F} \neq \emptyset$. On the other hand, given any chain $F_i$,
    $F_i \in \cc{F}$, the union $F = \bigcup_i\, F_i$ is clearly a
    inf-lattice filter such that $a \in F$. But $U \cap F = U \cap
    \bigcup_i \, F_i 
    = \bigcup_i (U \cap F_i) = \bigcup_i\, \emptyset =
    \emptyset$. Thus $F \in \cc{F}$. There exists then (by choice) a
    maximal element $P \in \cc{F}$. We show now that $P$ is
    $\jmath$-prime.

\vspace{1ex}

Given a inf-lattice filter $F \subset H$, and a element $a
\in H$, we denote \mbox{$( F, \; a ) = \{x \,|\, \exists \,b \in
F,\; b \wedge 
a \leq x\}$} the inf-lattice filter \emph{generated} by $F \cup \{a\}$.

\vspace{1ex}

Let $a_i \leq a \in \jmath(a)$ be a (finite non empty) cover, and [ (1)
  $a \in P$ ]. 
Assume that $\forall 
\,i\,\; a_i \notin P$. Then, $( P,\; a_i ) \cap U \neq \emptyset$.
Take $x_i \in U$, $x_i \in ( P,\; a_i )$, \mbox{$x_i \leq b_i
\wedge a_i$,} $b_i \in P$. It follows [ (2) $b_i \wedge a_i \in U$ ]. Let
$c = \bigwedge_i \, b_i$. \mbox{Then [ (3) $c \in P$ ]. } But $c \leq b_i\,$,
 thus $c \wedge a_i \leq b_i \wedge a_i$. It follows
from (2) that [ (4) $c \wedge a_i \in U$ ]. Since  $a_i \leq a \in
\jmath(a)$, we have $c \wedge a_i \leq c \wedge a \in \jmath(c \wedge
a)$ \mbox{(see \ref{gtopology} iii).} Thus, from (4) it follows [ (5) $c \wedge a \in
U$ ] (see \ref{2sheaves} ii). But from (1)
and (3), we have \mbox{[ (6) $c \wedge a \in P$ ].} Finally, (5) and
(6) contradict 
$U \cap P = \emptyset$. 
The case of the empty cover is easy. If $\emptyset \in \jmath(a)$, this
contradicts the fact that $a \notin U$.   
\end{proof} 

As a corollary of the previous two theorems we state

\begin{theorem} \label{spectralspace}
Let $(H,\, \jmath)$ be a site whose covers are \emph{finite
    families}. Consider the set $P_\jmath(H)$ of $\jmath$-prime
  inf-lattice filters $P \subset H$. Then, the sets $W_a = \{P \,|\, a
  \in P\}$ are compact and form a open base for a topology. The
  resulting topological 
  space $P_\jmath(H)$ 
  is sober, compact, and has a base of compact open sets. Its locale
  of open sets is (isomorphic to) the locale $I_\jmath(H)$ of
  $\jmath$-ideals of $H$ 
\qed\end{theorem}

Compact sober topological spaces with a base of compact opens are
known as \emph{spectral spaces} \cite{HO}, and first arisen as spaces of prime
ideals in ring theory.
\vspace{1ex}

All the results in this appendix apply to the following example:

\begin{example} \label{lattice}
Given any distributive lattice $V$, the \emph{finite suprema} form a
subcanonical Grothendieck topology (distributivity amounts to axiom
iii)), that we shall denote $\jmath_f$. The points $p$ of the site
$(V,\, \jmath_f)$ correspond with the prime filters $P \subset V$ of the
lattice, and a $\jmath_f$-ideal $U \subset V$ is just a lattice ideal
(notice that $0 \in U$ since $0$ is the empty supremun). 
The generated lattice ideal is given by:
$$
\hspace{6ex} \# S = \{x \,|\, \exists \;
  a_1,\, a_2,\; \ldots \; a_n \; \in S, \; n \geq 0,\; x \leq a_1\vee a_2\vee
  \; \ldots \; \vee a_n \}
$$ 
\qed\end{example}

\vspace{4ex}

{\bf  Eduardo J. Dubuc}

     Departamento de Matematicas, F.C.E. y N.,

     Universidad de Buenos Aires,

     Buenos Aires, Argentina.

\vspace{1ex}

{\bf  Yuri A. Poveda}

     Departamento de Matematicas,

     Universidad Nacional de Pereira,

     Pereira, Colombia.


\begin{thebibliography}{99}

%\bibitem{G2} Artin M, Grothendieck A, Verdier J.,  \textsl{SGA 4 ,
%(1963-64)}, Springer Lecture Notes in Mathematics 269 (1972).

%\bibitem{G3} Artin M, Grothendieck A, Verdier J.,  \textsl{SGA 4 ,
%(1963-64)}, Springer Lecture Notes in Mathematics 270 (1972).

\bibitem{COM} Cignoli R., D'Ottaviano I., Mundici
  D., \textsl{Algebraic Foundations of Many-valued Reasoning}, Trends
  in Logic Vol 7, Kluwer Academic Puplishers (2000).

\bibitem{CDM} Cignoli R., Dubuc E. J., Mundici
  D., \textsl{Extending Stone duality to multisets and locally finite
    MV-algebras}, Journal of Pure and Applied Algebra 189 (2004).

\bibitem{C} Coste M., \textsl{Localization, Spectra and Sheaf Representation}, Springer Lecture Notes in Mathematics 753 (1977). 
 
\bibitem{H} Hartshorne R., \textsl{Algebraic Geometry}, Springer
  Verlag, New York (1977).

\bibitem{HO} Hochster M., \textsl{Prime ideal structure in commutative
  rings}, Trans. Amer. Math. Soc. 142, 43-60 (1969).

\bibitem{J} Johnstone P. T., \textsl{Stone spaces}, Cambridge
  University Press (1982).

\bibitem{JT} Joyal A., Tierney M., \textsl{An extension of the Galois
  Theory of Grothendieck}, Memoirs of the American Mathematical
  Society, Vol. 151, (1984). 

\bibitem{M} Mac Lane S., \textsl{Categories for the working
  mathematician}, 2nd ed., Springer Verlag, (1998).

\bibitem{MM} Mac Lane S., Moerdijk I., \textsl{Sheaves in Geometry and
  Logic}, Springer Verlag, (1992).

\bibitem{NPM} Novak V., Perfilieva I., Mockor J., \textsl{Mathematical Principles of Fuzzy Logic}, Kluwer Academic Publishers, Boston (1999).

\end{thebibliography}
\end{document}